\documentclass{article}

\usepackage[colorlinks,bookmarksopen,bookmarksnumbered,citecolor=red,urlcolor=red]{hyperref}
\usepackage[a4paper,body={18cm,26cm},top=2cm]{geometry}

\usepackage{amsthm,amsmath,amssymb,stmaryrd,type1cm}

\usepackage{pgfplots}
\pgfplotsset{compat=1.3}
\usepackage{pgfplotstable}
\usepackage[caption=false]{subfig}
\usepackage{booktabs}

\newtheorem{remark}{Remark}
\newcommand{\V}[1]{\underline{#1}}                     
\newcommand{\M}[1]{\underline{\underline{#1}}}
\newcommand{\Vt}[1]{\underline{\tilde{#1}}}                     
\newcommand{\Mt}[1]{\underline{\underline{\tilde{#1}}}}

\newcommand{\dofP}{\mathbb{V}}
\newcommand{\shapeP}{\V{\mathfrak{V}}}
\newcommand{\forceP}{\mathbb{F}}

\newcommand{\Domain}{\Omega}
\newcommand{\domain}{\omega}
\newcommand{\Inte}{\Gamma}
\newcommand{\inte}{\gamma}
\newcommand{\DomI}{\Domain_I}
\newcommand{\DomB}{\Domain_B}

\newcommand{\domI}{\domain_I}
\newcommand{\domB}{\domain_B}
\newcommand{\domC}{\domain_C}

\newcommand{\IntC}{\Inte_C}
\newcommand{\intI}{\inte_I}
\newcommand{\intC}{\inte_C}

\def\traitemr #1,#2\endd{\ifx\end #1\else\ref{#1}\ \traitemr #2\endd\fi}

\graphicspath{{figures/}}

\title{Nonintrusive coupling of 3D and 2D laminated composite models based on finite element 3D recovery}

\author{Guillaume Guguin$^{1}$, Olivier Allix$^{1}$,\\ Pierre Gosselet$^{1}$, St\'{e}phane Guinard$^{2}$\\[10pt]
 \small (1)LMT-Cachan, ENS-Cachan/CNRS/Pres UniverSud Paris,\\ \small61 avenue du Pr\'esident Wilson, 94235 Cachan, France\\\small
(2)  EADS-IW, Toulouse, France }

\date{Jan. 2014, IJNME \\ doi:10.1002/nme.4630 }

\begin{document}
\maketitle

\begin{abstract}
In order to simulate the mechanical behavior of large structures assembled from thin composite panels, we propose a coupling technique which substitutes local 3D models for the global plate model in the critical zones where plate modeling is inadequate.
The transition from 3D to 2D is based on stress and displacement distributions associated with Saint-Venant problems which are precalculated automatically for a simple 3D cell. The hybrid plate/3D model is obtained after convergence of a series of iterations between a global plate model of the structure and localized 3D models of the critical zones. This technique is nonintrusive because the global calculations can be carried out using commercial software. Evaluation tests show that convergence is fast and that the resulting hybrid model is very close to a full 3D model.

\textbf{keywords: }{Plate-3D coupling; mixed dimensionality; nonintrusive coupling; laminated composites.}

\end{abstract}

\section{Introduction}

This paper deals with nonintrusive techniques for the simulation of structures assembled from large flat composite panels. In many such structures, one dimension is much smaller than the other two, so they can be modeled as plates. However, plate theories alone cannot take into account the 3D stress concentration states which occur near edges, geometric discontinuities or zones subjected to concentrated loads. These phenomena are important in the case of laminated composites because they are often responsible for the initiation and propagation of damage, e.g. in the case of  delamination  \cite{garg88}. Therefore, a natural approach to dealing with large composite structures is to attempt, whenever necessary, to couple plate or shell models with 3D models \cite{Duster20073524}.

The main topic of this paper is the development and evaluation of a type a coupling between a plate model and a 3D model which could be applied directly within the nonintrusive framework proposed in \cite{GENDRE.2009.1,GENDRE.2011.1} {as well as in \cite{bettinotti2013coupling} for dynamics}. {Essentially, the nonintrusive coupling technique consists, starting from an initial global calculation, in alternating between local calculations with prescribed displacements at the boundaries and global calculations which include corrective loads in order to reduce the imbalance between the models.}
Because the nonintrusive coupling technique used in the paper was designed for models with the same dimensionality, the main problem addressed here is the recovery of the 3D quantities associated with the 2D model.

Overall, the question of the link between a plate model and its 3D counterpart is of such practical and theoretical importance that it has given rise to many publications. This introduction focuses more specifically on the references of interest in the literature that we tried to incorporate into our formulation.

The  quality of a plate model is usually assessed  by comparison with a corresponding 3D reference model. This has led to a first group of works focusing on the development of \textit{a posteriori} error estimators, initially neglecting edge effects \cite{koi70,dan70,sim72,koi72b,lad76}, then taking them into account \cite{lad80}. Thus, by the end of the 1980's, the performance of the classical Kirchhoff-Love and Reissner-Mindlin plate and shell theories in 3D elasticity was relatively well-understood. Another group of works concerned the technique of asymptotic developments with the thickness as the small parameter \cite{fri61,gol66,cia79,des86a,des86b,cai87,cia90}. In addition, using the same approach as Reissner \cite{reissner45,reissner86} in statics and Mindlin in dynamics \cite{mindlin51}, a shear factor was often introduced \cite{whitney73} in order to correct what would have been deduced from a pure kinematic assumption by taking into account some \textit{a priori} information on the three-dimensional stress state. This mixed nature of the link between plate theory and 3D theory is the reason why the construction of improved theories  \cite{Lur64,tou91,lev80,red84a,carrera02} and associated finite elements \cite{carrera04,auricchio98} relies heavily on mixed formulations.

Regardless of the method, the assessment of a plate model involves the comparison between a 3D ``Saint-Venant'' or long-wavelength solution \cite{ladeveze83,gregory85} and the solution obtained from plate theory. A key aspect of the comparison, especially in the case of composites, concerns the distribution of the stress field  \cite{pagano78,davi96}. Therefore, a basic technique which is used, sometimes implicitly, when estimating the quality of plate theories is the recovery of equilibrated 3D solutions by integrating the equilibrium equations through the thickness \cite{rolfes98,o.allix_plate_2009}. The application of this technique to finite element calculations requires the recovery of 2D fields which satisfy the 2D equilibrium equations exactly. The use of such a stress recovery technique \cite{noor90,daghia_hybrid_2008} makes the implementation within a commercial program difficult. A comparative analysis of various models and techniques to evaluate \textit{a priori} or \textit{a posteriori} out-of-plane stresses was presented in \cite{carrera00}. An improved hybrid post-processing procedure was proposed in \cite{cen02} along with a new displacement-based element. This procedure, which works element-wise and, thus, is easier to implement, was extended to first-order theories in \cite{kim07} to recover both displacements and stresses.

In order for the 3D approximation derived from plate theory to be valid everywhere, it should be completed by adding an approximation of the edge effects. A key issue is that the problem deduced from the plate theory residual at the edges must lead to a localized problem \cite{lad80,lad87a,o.allix_plate_2009}. Otherwise, even using what one calls refined theories, the plate's interior solution would be meaningless for the orders of magnitude considered. This may be one of the reasons why the Reissner-Mindlin theory is the dominant approach today, even for composite structures. Another reason is, of course, that because it involves only first-order derivatives it is much easier to build finite elements based on this theory than, for example, on Kirchhoff-Love's. However, the Reissner-Mindlin theory includes an artificial edge effect which, if not handled properly, could pollute the recovery of the 3D solution near the edges \cite{rossle_corner_2011}.

In this context, our approach to nonintrusive coupling must comply with several constraints. The first constraint is that our method involves two overlapping models, a 2D (plate or shell) model which is defined over the whole domain to be analyzed and remains unchanged throughout the analysis, and a refined 3D model in which structural details and complex nonlinearities such as contact or delamination can be taken into account. The second constraint is that we need to use results obtained from commercial finite element codes which, usually, are based on the Reissner-Mindlin model and do not provide recovery techniques. The implementation of such recovery techniques would require rather complex manipulations with simultaneous access to data from several elements, which would annihilate the very purpose of nonintrusive techniques. Moreover, the plate model is usually less than optimal compared to its 3D counterpart in that the plate quantities which would be deduced from the calculation of the full 3D model would be somewhat different from those obtained using the plate model, even far away from the edges. That is the reason why \cite{noor90,yu03} and others have used an asymptotic formulation up to the second order to define a quasi-optimal version of Reissner-Mindlin's composite plate theory by means of an optimization procedure. Then, based on these accurate plate quantities, a recovery technique can be used to obtain the 3D displacement, strain and stress fields as functions of 2D variables.

{
The paper is organized as follows:
Section~\ref{sec:recovery} shows how relevant 3D quantities can be generated for any stacking sequence in order to recover 3D stress and displacement bases. Then, in Section~\ref{sec:models}, these bases are used to define hybrid models.
In Section~\ref{sec:study}, a first analysis of the results enables us to identify the most favorable coupling technique, which is tested on a more representative example in Section~\ref{sec:assess}.
}

\section{Numerical recovery of the Saint-Venant stress and warping}\label{sec:recovery}
In this section we investigate the association of 3D quantities (traction and warping) with composite plate calculations.

\subsection{Notations and equations}

Let us consider a classical 3D problem defined in domain $\Domain$.
Let $\V{f}$ be the volume force, $\V{g}$ the traction prescribed over $\partial_N\Domain$ and $\V{u}_D$ the prescribed displacement over the complementary part $\partial_D\Domain$. {The unit vectors of an orthogonal spatial basis are denoted by $(e_{x},e_{y},e_{z})$}. The { linear elastic behavior} is given by Hooke's tensor $\mathbb{H}$.
The symmetric part of the displacement gradient and the Cauchy stress tensor are designated respectively by $\M{\epsilon}(\V{u})$ and $\M{\sigma}$.

We introduce the affine space of the admissible displacements $\mathcal{U}(\Domain)$:
\begin{equation*}
\mathcal{U}(\Domain) = \left\{ \V{u}\in H^1(\Domain),\ \V{u}=\V{u}_D \text{ on } \partial_D\Domain\right\}
\end{equation*}
whose associated vector space is $\mathcal{U}_0(\Domain)$.

The 3D problem is:
\begin{equation}
\left\{\begin{aligned}
&\text{Find }\V{u}\in \mathcal{U}(\Domain) \text{ such that, }\forall \V{u}^*\in\mathcal{U}_0(\Domain),  \\
&\int_\Domain \M{\sigma}:\M{\epsilon}(\V{u}^*)dx = \int_\Domain \V{f}\cdot\V{u}^*dx +  \int_{\partial_N\Domain} \V{g}\cdot\V{u}^*dS \\
&\M{\sigma}=\mathbb{H}:\M{\epsilon}(\V{u})
\end{aligned}\right.
\end{equation}

Regarding the plate model, we have $\Domain=\domain\times [-h,h]$.

Using the notation $\tilde{\ }$ to designate in-plane quantities ($x$ and $y$ coordinates), the first assumption concerns the shape of the plate's displacement:
\begin{equation}
\V{v} = \Vt{v}+v_z \V{e}_z +z\Vt{\theta}\wedge \V{e}_z
\end{equation}
where $\Vt{v}$ and $v_z$ are respectively the {in-plane and out-of-plane displacements of the mid-surface} and $\Vt{\theta}$ denotes the {rotation of the section}, all of which are functions defined in $\domain$. {Operator $\wedge$ is the cross product defined by:
\[\begin{pmatrix}
x_{1}\\
x_{2}\\
x_{3}
\end{pmatrix}\wedge\begin{pmatrix}
y_{1}\\
y_{2}\\
y_{3}
\end{pmatrix}=\begin{pmatrix}
x_{2}y_{3}-x_{3}y_{2}\\
x_{3}y_{1}-x_{1}y_{3}\\
x_{1}y_{2}-x_{2}y_{1}
\end{pmatrix}
\]
}
This leads to the introduction of the classical plate strains:
\begin{equation}
\M{\epsilon}(\V{v}^P)=
 \begin{pmatrix}
\Mt{\epsilon}& \Vt{\epsilon}_z \\ \text{sym} & v_{z,z}
\end{pmatrix} \text{ with } \left\{\begin{aligned}
&\Mt{\epsilon} = \Mt{\gamma}(\Vt{v}) + z \Mt{\chi}(\Vt{\theta}) \\
&\Vt{\epsilon}_z = \frac{\tilde{\V{\nabla}}(v_z)+\Vt{\theta}\wedge\V{e}_z}{2}
\end{aligned}\right.
\end{equation}
\begin{equation*}
\Mt{\gamma}(\Vt{v}) = \begin{pmatrix}
v_{x,x} & \frac{ v_{x,y} + v_{y,x} }{2} \\ \text{sym} & v_{y,y}
\end{pmatrix}, \qquad \Mt{\chi}(\Vt{\theta}) = \begin{pmatrix}
\theta_{y,x} & \frac{\theta_{y,y}-\theta_{x,x}}{2} \\ \text{sym} & -\theta_{x,y}
\end{pmatrix}
\end{equation*}

Using notation $\langle \rangle$ for the integration through the thickness, the plate's stress tensors are:
\begin{equation}\label{eq:stressplate0}
\Mt{N} = \langle \Mt{\sigma} \rangle\qquad
\Mt{M} = \langle z\Mt{\sigma} \rangle \qquad
 \Vt{Q} = \langle \Vt{\sigma}_z \rangle
\end{equation}

The second assumption of the plate model is that the peeling stress $\sigma_{zz}$ is negligible compared to the other components of the stress tensor. This assumption {simplifies} the expression of the {plate's} behavior as a function of the plane stress $\tilde{\mathbb{H}}_{ps}$ {and out-of-plane stress} $\tilde{\mathbb{B}}$. (In 3D, $\tilde{\V{\sigma}}_z=\tilde{\mathbb{B}}{\Vt{\epsilon}}_z$.)

For a plate in equilibrium, the constitutive relations are:
\begin{equation}
\begin{aligned}
\Mt{N}&=\langle\tilde{\mathbb{H}}_{ps}\rangle : \Mt{\gamma}(\Vt{v}) \\
\Mt{M}& = \langle z^2\tilde{\mathbb{H}}_{ps}\rangle : \Mt{\chi}(\Vt{\theta})\\
\Vt{Q}& = \kappa \langle\tilde{\mathbb{B}}\rangle :\Vt{\epsilon}_z(v_z,\Vt{\theta})
\end{aligned}
\end{equation}
where $\kappa$ is a shear correction factor. (In this paper, we use the correction factor from \cite{batoz_moderation_1990}.)

The admissible plate displacement field is:
\begin{equation}
\mathcal{U}^P(\domain)=\left\{(\Vt{v},\Vt{\theta},v_z)\in H^1(\domain)^5, \text{which satisfies Dirichlet's BCs}\right\}
\end{equation}
and the plate problem becomes\footnote{for the sake of simplicity, we did not develop the right-hand-side in terms of plate components}:
\begin{equation}
\left\{
\begin{aligned}
&\text{Find }(\Vt{v},\Vt{\theta},v_z)\in \mathcal{U}^P(\domain) \text{ such that }\forall (\Vt{v}^*,\Vt{\theta}^*,v_z^*)\in\mathcal{U}_0^P(\domain)  \\
&\!\!\int_\domain\left(\Mt{N}:\Mt{\gamma}(\Vt{v}^*) + \Mt{M}:\Mt{\chi}(\Vt{\theta}^*)+ \Vt{Q}\cdot\Vt{\epsilon}_z(v_z^*,\Vt{\theta}^*)\right)dx = \int_\Domain \V{f}\cdot\V{v}^* dx +  \int_{\partial_N\Domain} \V{g}\cdot\V{v}^*dS  \\
&\Mt{N} = \langle\tilde{\mathbb{H}}_{ps}\rangle : \Mt{\gamma}(\Vt{v}),\quad \Mt{M} = \langle z^2\tilde{\mathbb{H}}_{ps}\rangle : \Mt{\chi}(\Vt{\theta}),\quad \Vt{Q} =\kappa  \langle\tilde{\mathbb{B}}\rangle :\Vt{\epsilon}_z(v_z,\Vt{\theta})
\end{aligned}\right.
\end{equation}

\begin{remark}
The Reissner-Mindlin kinematic assumption was chosen because it is the most commonly used in commercial software. This model provides an estimation of the transverse shear forces, but, unfortunately, it leads to perturbations of the generalized forces near the edges \cite{rossle_corner_2011}.
\end{remark}

\begin{remark}[Notations]
From now on, we will use the following notations:

for the kinematics:
\begin{equation*}
\left.\begin{aligned}
\text{plate DOFs } \dofP &=\begin{pmatrix} v_x& v_y & \theta_x&\theta_y& v_z \end{pmatrix}^T \\
\text{plate kinematics }\shapeP&=\begin{pmatrix} \V{e}_x & \V{e}_y & -z\V{e}_y &z\V{e}_x & \V{e}_z\end{pmatrix} \end{aligned}\right|
\text{ so that }\V{v} = \sum_{k=1}^5 \shapeP_k \dofP_k
\end{equation*}
and for the generalized forces:
\begin{equation}\label{eq:stressplate}
\begin{aligned}
\forceP &= \begin{pmatrix} N_{xx} & N_{xy} & M_{xx} & M_{xy} & Q_x & N_{yy} & M_{yy} & Q_y\end{pmatrix}^T \\ \text{so that }\quad
\forceP_{1:5} &= \langle \M{\sigma}:(\shapeP_{1:5}\otimes \V{e}_x ) \rangle \\
\forceP_{6:8} &= \langle \M{\sigma}:(\shapeP_{2,4,5}\otimes \V{e}_y ) \rangle
\end{aligned}
\end{equation}
Although these notations seem biased toward the $(\V{e}_x,\V{e}_y)$ coordinate system (especially toward direction $\V{e}_x$ since the components related to this normal direction are listed first), the following calculations restore the normal-independent property of the method.

In order to deal with interfaces which are not aligned with the axes, we introduce the following notations: for an interface whose local coordinate system is $(\V{n},\V{t},\V{e}_z)$ ($\V{n}$ being the normal vector and $\V{t}$ the tangent vector), $\phi$ denotes the angle such that $\cos(\phi)=\V{n}\cdot\V{e}_x=n_x$. The kinematic basis of the plate in the local coordinate system is $\shapeP^\phi=(\V{n},\V{t},-z\V{n},z\V{t},\V{e}_z)$. The associated DOFs are $\dofP^\phi=(v_n,v_t,\theta_n,\theta_t,v_z)$ and the stress components are $\forceP^\phi= \begin{pmatrix} N_{nn} & N_{nt} & M_{nn} & M_{nt} & Q_n & N_{tt} & M_{tt} & Q_t\end{pmatrix}^T$.
\end{remark}

\subsection{Recovery of the Saint-Venant traction fields}
Since plate solutions are relevant only for problems with long characteristic lengths of variation, we try to associate one 3D Saint-Venant stress field $(\M{\tau}_i)$ with each plate stress component.

\begin{figure}[htbp]\centering
   \begin{minipage}[c]{.49\linewidth}\centering
   \includegraphics[width=.99\linewidth]{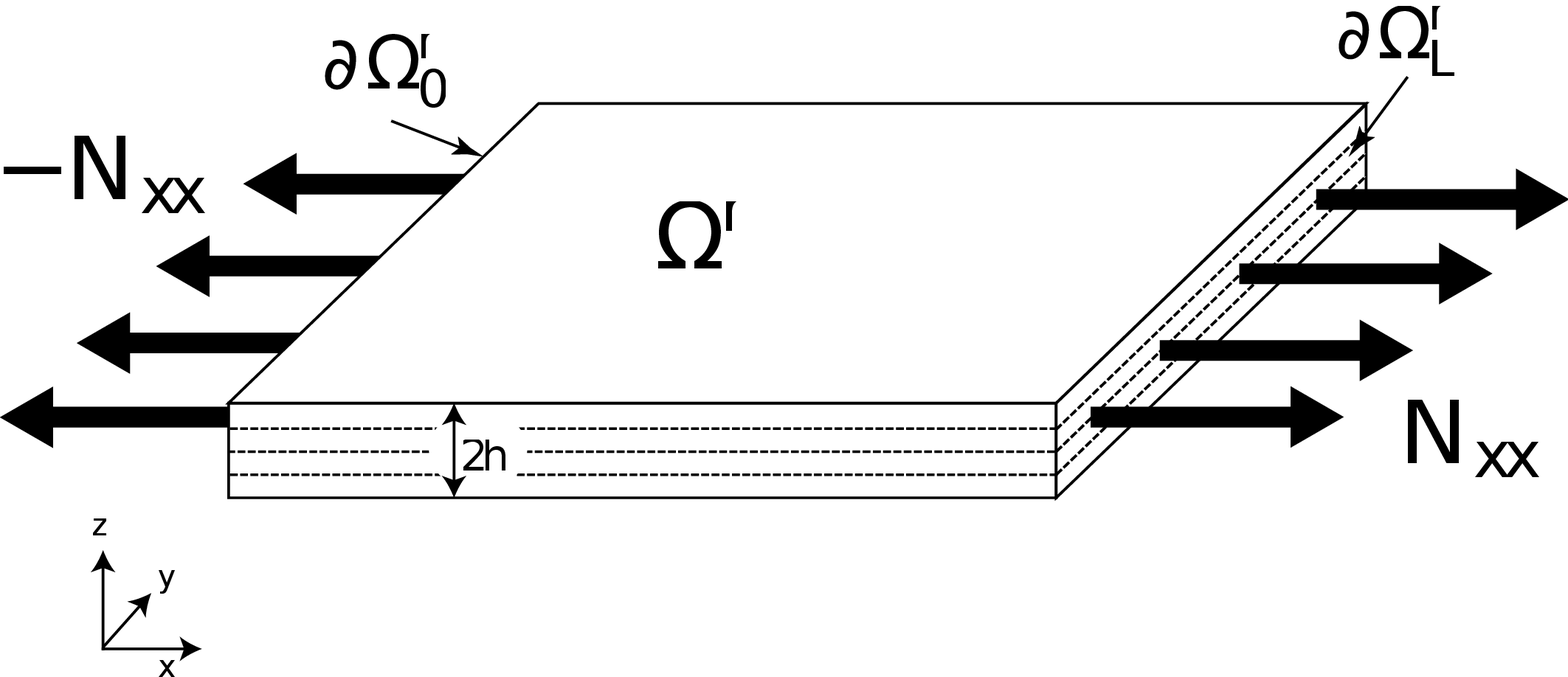}
   \caption{Cell problem, loading $ N_{x} $}\label{fig:PbRelev1}
   \end{minipage}
      \begin{minipage}[c]{.49\linewidth}\centering
   \includegraphics[width=.99\linewidth]{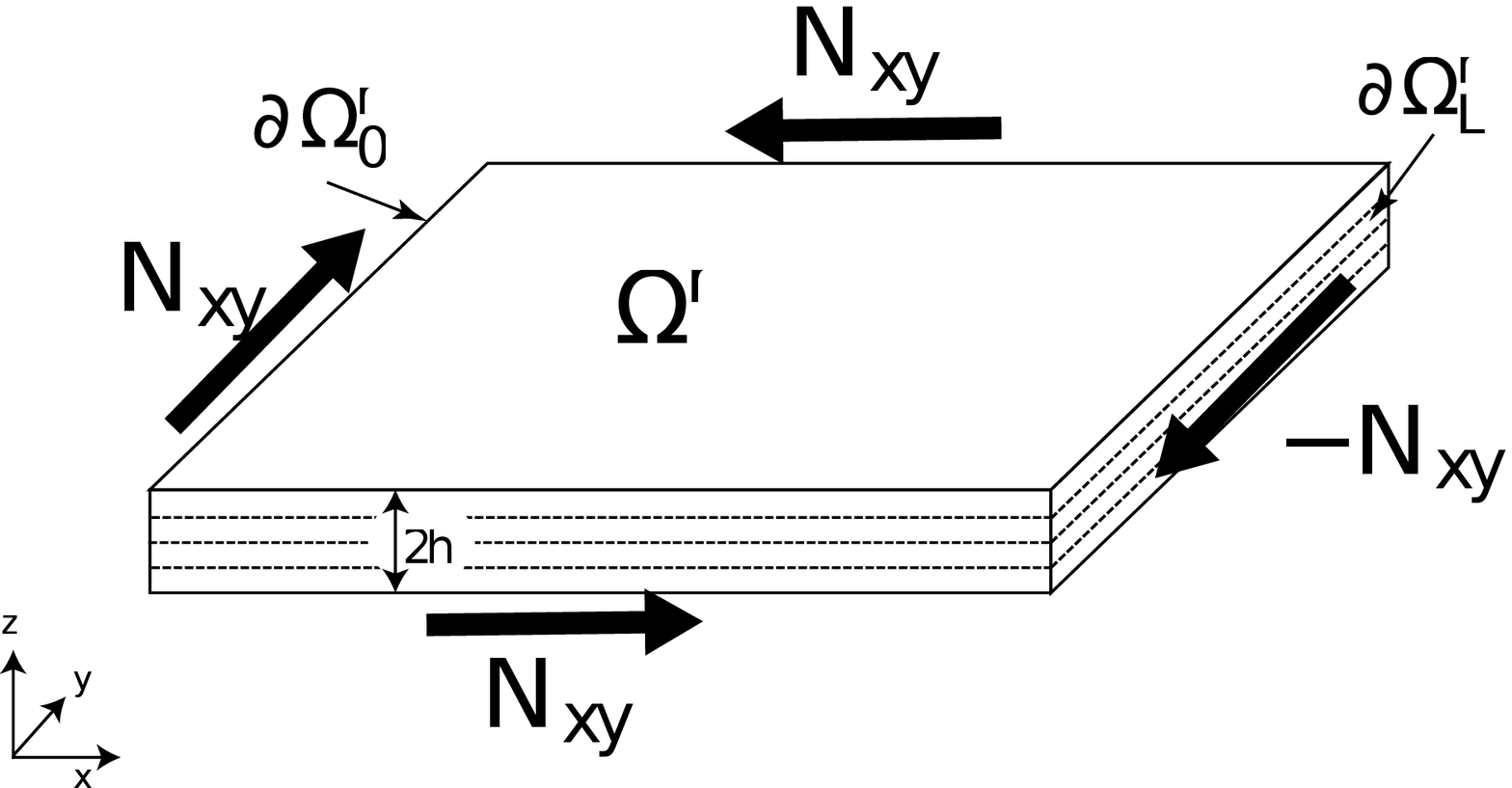}
      \caption{Cell problem, loading  $ N_{y} $}\label{fig:PbRelev2}
   \end{minipage}
      \begin{minipage}[c]{.49\linewidth}\centering
   \includegraphics[width=.99\linewidth]{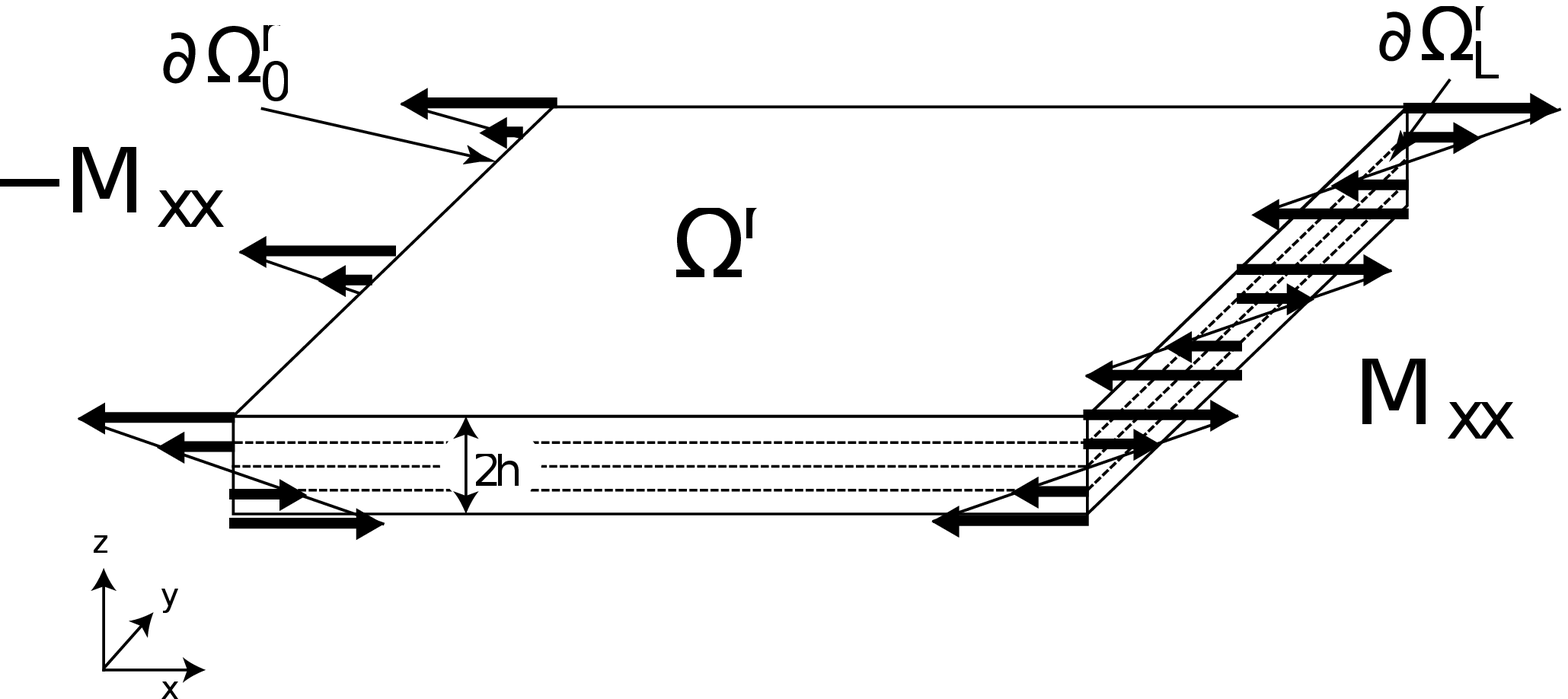}
      \caption{Cell problem, loading  $ M_{x} $}\label{fig:PbRelev3}
   \end{minipage}
      \begin{minipage}[c]{.49\linewidth}\centering
   \includegraphics[width=.99\linewidth]{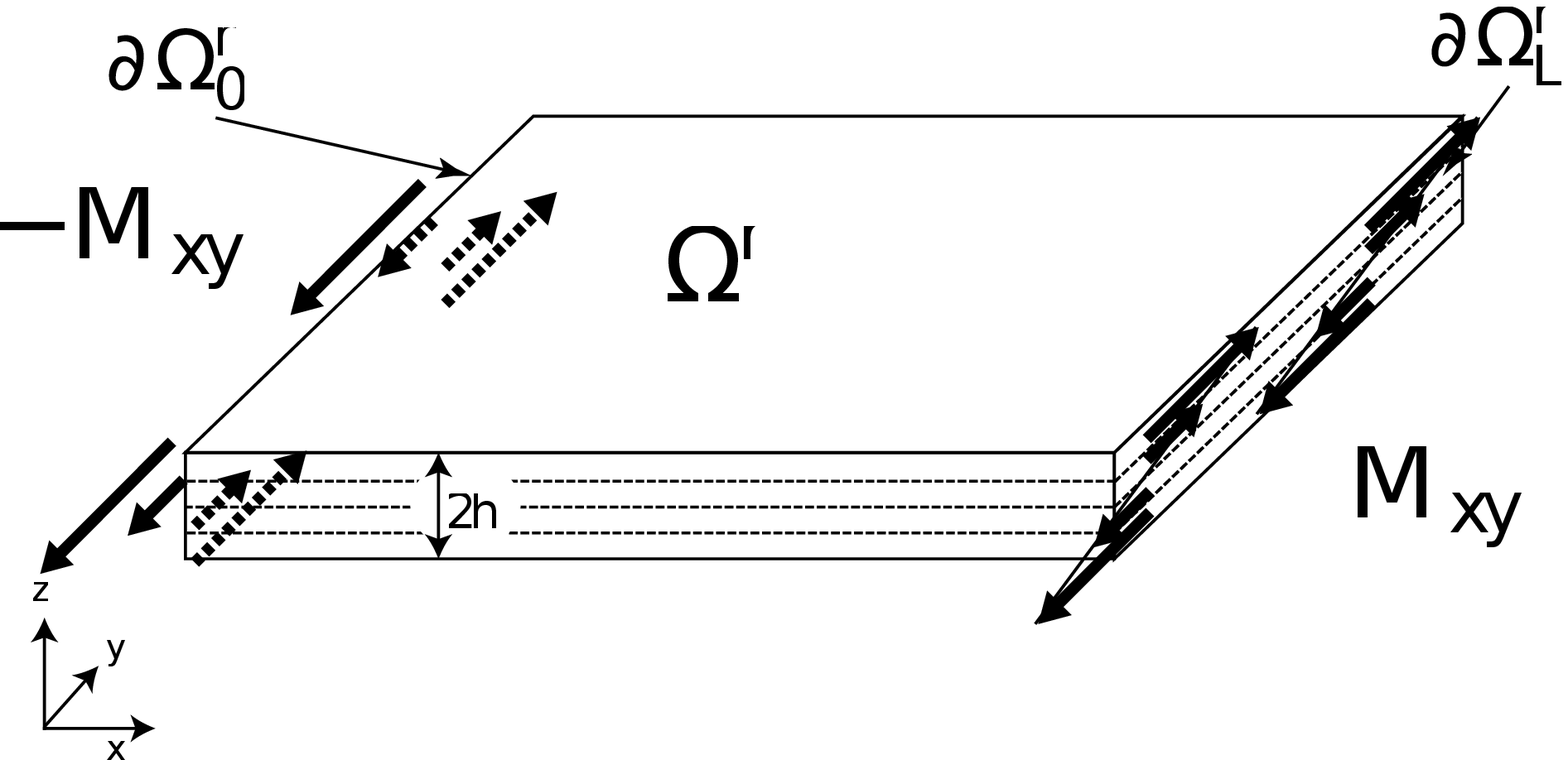}
     \caption{Cell problem, loading  $ M_{y} $}\label{fig:PbRelev4}
   \end{minipage}
      \begin{minipage}[c]{.49\linewidth}\centering
   \includegraphics[width=.99\linewidth]{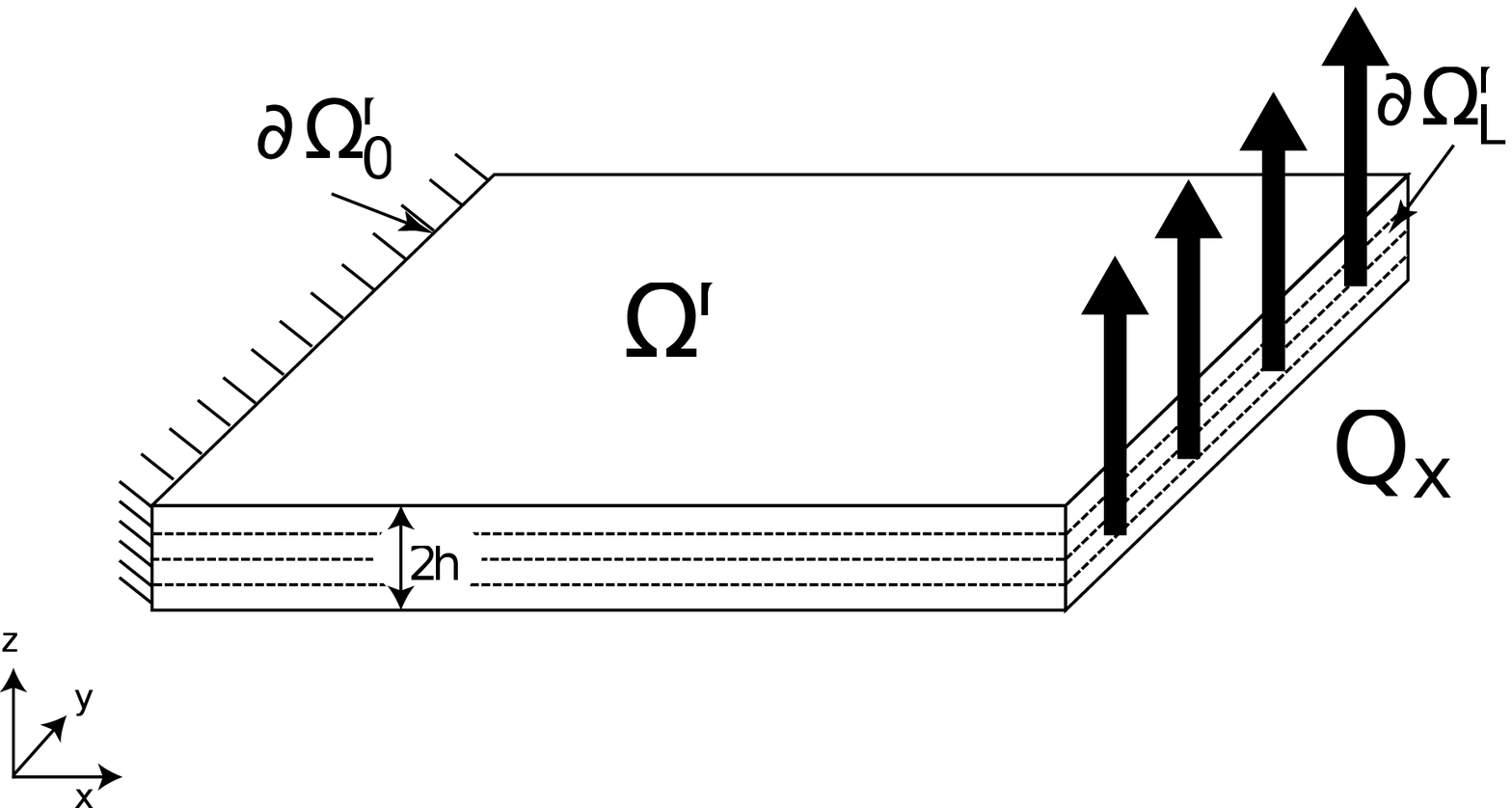}
     \caption{Cell problem, loading  $ Q $}\label{fig:PbRelev5}
   \end{minipage}
\begin{minipage}[c]{.49\linewidth}\centering
   \includegraphics[width=.99\linewidth]{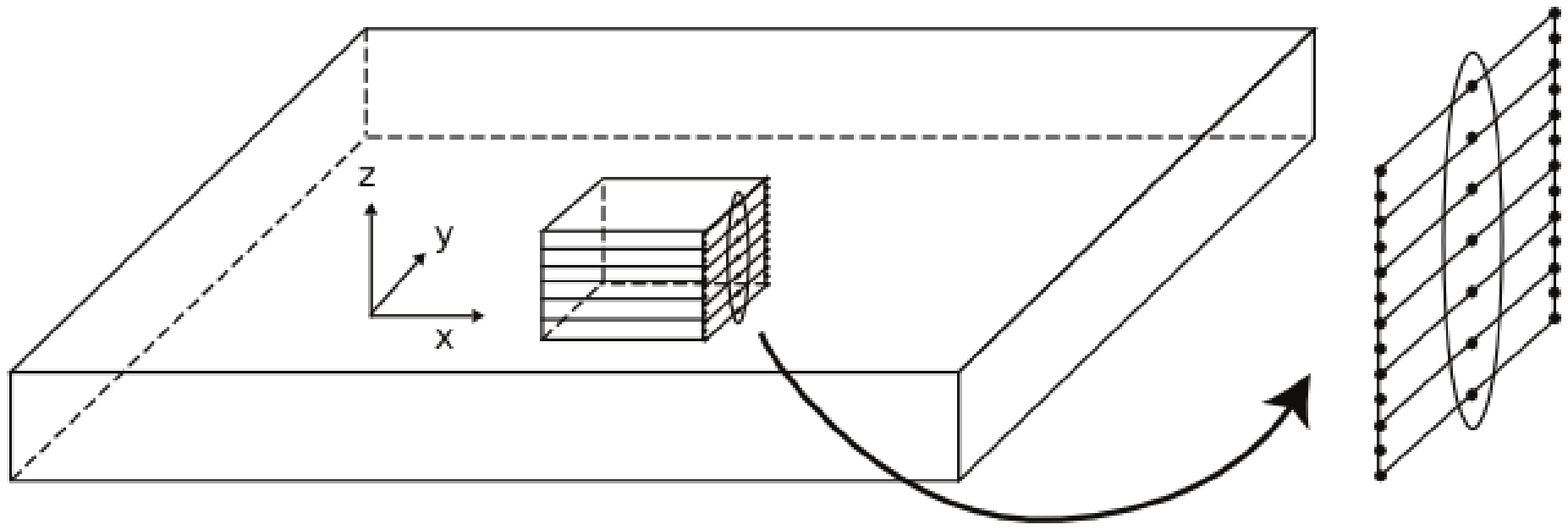}
     \caption{Extraction of the stress profiles and nodal forces}\label{fig:Extract}
   \end{minipage}
\end{figure}

{
In order to do that, we consider the simplest 3D problems whose solutions are Saint-Venant fields, which consist in a sufficiently large square plate with the same material, the same stacking sequence and the same finite elements as the local model, subjected to 8 components of generalized forces (3 in tension, 3 in bending and 2 in shearing). This leads to 8 basic problems with well-chosen regular loads, each activating one generalized force as defined in Table~\ref{tab:CLrelev} and illustrated in Figures~\ref{fig:PbRelev1}-\ref{fig:PbRelev5}.
We call these problems ``cell problems'' because of the analogy of the method with micro-macro approaches. Here, we are less interested in the average response than in the distribution of the 3D mechanical quantities through the thickness in the center of the plate, which represents a numerical ``3D recovery'' of the plate quantities.}

\begin{table}[ht]\centering
\begin{tabular}{|c|c|c|c|}\hline
j & loading & main {load} & figure \\\hline
1 & $\sigma_{xx}  = -1  \text{ over }\partial\Domain_{0},\  \sigma_{xx}  = 1  \text{ over }\partial\Domain_{L} $ & $N_{xx}$ & \ref{fig:PbRelev1}\\\hline
2  & $\begin{aligned}
&\sigma_{xy}  = -1  \quad\text{over}\,\partial\Domain_{L}  & \sigma_{xy}=1 & \quad\text{over}\,\partial\Domain_{0} \\
&\sigma_{xy} = -1  \quad\text{over}\,\partial\Domain_{a}  & \sigma_{xy} = 1& \quad\text{over}\,\partial\Domain_{-a}
\end{aligned}$ & $N_{xy}$ & \ref{fig:PbRelev2}\\\hline
3 & $\sigma_{xx} = z \text{ over }\partial\Domain_{0},\  \sigma_{xx} = z  \text{ over }\partial\Domain_{0}$ & $M_{xx}$& \ref{fig:PbRelev3} \\\hline
4 & $\begin{aligned}
&\sigma_{xy} = -z  \quad\text{over}\,\partial\Domain_{0} && \sigma_{xy} = z  \quad\text{over}\,\partial\Domain_{L}\\
&\sigma_{xy} = z  \quad\text{over}\,\partial\Domain_{-a} && \sigma_{xy} = -z  \quad\text{over}\,\partial\Domain_{a}
\end{aligned}$ & $M_{xy}$ &\ref{fig:PbRelev4} \\\hline
5 & $\sigma_{xz}   = 1 \quad\text{over}\,\partial\Domain_{L}$ & $Q_x$  &\ref{fig:PbRelev5}\\\hline
6& same as Problem 1, but rotated $90^\circ$ & $N_{yy}$&\\\hline
7& same as Problem 3, but rotated $90^\circ$ &$M_{yy}$& \\\hline
8& same as Problem 5, but rotated $90^\circ$ & $Q_{y}$&\\\hline
\end{tabular}\caption{Construction of the Saint-Venant solutions}\label{tab:CLrelev}
\end{table}

From each calculation $j\in\llbracket 1,8\rrbracket$, one can extract the distribution of the 3D stress tensor $(\M{\sigma}_j)$ through the thickness in the center of the plate, \emph{i.e.} sufficiently far away from any edge effect. This is illustrated in Figure~\ref{fig:Extract} after FE discretization. From this stress tensor, using Equation~\eqref{eq:stressplate}, one can post-process the generalized forces $(\forceP_j)$ of the plate.
These test cases were designed so that $\forceP^j_j$ is dominant in the $j^{th}$ experiment (although coupling among generalized forces exists naturally, if only because of the conservation of moments which involves shear forces).

Because of the linearity of the problem, any Saint-Venant stress state is a combination of the extracted $(\M{\sigma}_j)$. The $(\M{\tau}_i)$ are constructed such that the contributions of the generalized forces are uncoupled:
\begin{equation}
\sum_{i=1}^8 \forceP^i_j \M{\tau}_i = \M{\sigma}_j \quad\Longrightarrow\quad \begin{bmatrix}
\M{\tau}_1 &\hdots & \M{\tau}_8
\end{bmatrix} = \begin{bmatrix}
\M{\sigma}_1 &\hdots &\M{\sigma}_8
\end{bmatrix} \left[\forceP^i_j\right]^{-1}
\end{equation}
where brackets denote matrices. Then, the $(\M{\tau}_i)$ are obtained by inversion of the $8\times 8$ matrix $\left[\forceP^i_j\right]$.

In other words, the $(\M{\tau}_i)$ belong to the space of the Saint-Venant stresses and satisfy the orthogonality relation:
\begin{equation}
\begin{aligned}
1 &= \langle \M{\tau}_i:(\shapeP_k\otimes \V{e}_v ) \rangle\ \begin{array}{l}\text{for } v=x \text{ and } i=k\in\llbracket 1,5\rrbracket\\ \text{or for } v=y\text{ and } (i,k)\in\{(6,1),(7,3),(8,5)\}\end{array} \\
0 &= \langle \M{\tau}_i:(\shapeP_k\otimes \V{e}_v ) \rangle\ \text{in all other cases }
\end{aligned}
\end{equation}

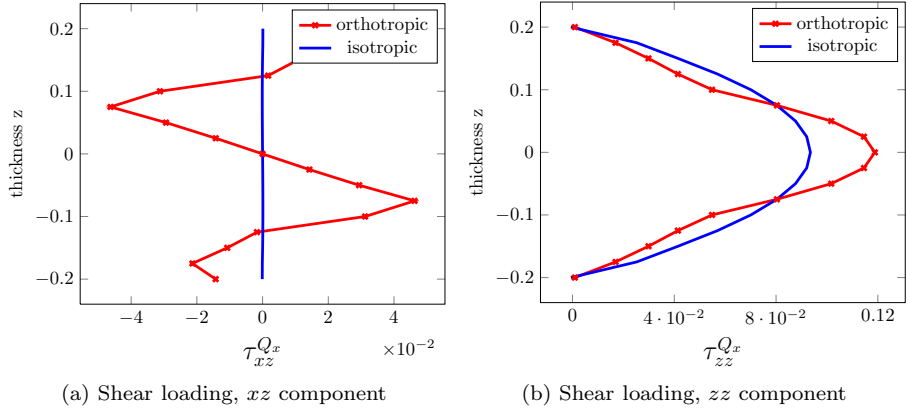
\begin{figure}[ht]\centering
\subfloat[][Shear loading, $xz$ component]{
   \begin{tikzpicture}[scale=0.7]
\begin{axis}[
		xlabel={\Large $\tau^{Q_x}_{xz}$},
		ylabel={thickness z},
		ylabel near ticks,tick scale binop=\times,
		scaled y ticks=true,]
		\addplot[color=red,mark=x,line width=1.5] table[x index=1,y index=0]{figures/relev/ortho_tau_QX.txt};
		\addlegendentry{orthotropic}
		\addplot[color=blue,line width=1.5] table[x index=1,y index=0]{figures/relev/iso_tau_QX.txt};
		\addlegendentry{isotropic}
\end{axis}
\end{tikzpicture}
\label{fig.tau.FX}
}
\subfloat[Shear loading, $zz$ component]{
\begin{tikzpicture}[scale=0.7]
\begin{axis}[
		xlabel={\Large $\tau^{Q_x}_{zz}$},xtick={0,0.04,0.08,0.12},
		ylabel={thickness z},
		ylabel near ticks,tick scale binop=\times,
		scaled y ticks=true,]
		\addplot[color=red,mark=x,line width=1.5] table[x index=3,y index=0]{figures/relev/ortho_tau_QX.txt};
		\addlegendentry{orthotropic}
		\addplot[color=blue,line width=1.5] table[x index=3,y index=0]{figures/relev/iso_tau_QX.txt};
		\addlegendentry{isotropic}
\end{axis}
\end{tikzpicture}
\label{fig.tau.FZ}
}
\caption[]{Distribution of the nodal forces through the thickness for isotropic and orthotropic materials}\label{fig.tau}
\end{figure}

Figure~\ref{fig.tau} shows typical distributions of the nodal forces through the thickness, both for an isotropic plate and for 4 orthotropic plies. (A complete description of the material will be presented in Section~\ref{sec:study}.)
In the isotropic case, we recover the analytical functions of \cite{mccune_mixed-dimensional_2000}.

\begin{remark}
The $(\M{\tau}_i)$ were designed for loads which are aligned with axes $(\V{e}_x,\V{e}_y)$. If necessary, one can substitute a tensor family $(\M{\tau}^\phi_i)$ {in a rotated direction}: given a local coordinate system $(\V{n},\V{t})$, let $\V{n}=n_x\V{e}_x+n_y\V{e}_y$ so that $\V{t}=-n_y\V{e}_x+n_x\V{e}_y$. Using linearity, the Saint-Venant solution for a load which is aligned with the local frame is:
\begin{equation}
\begin{aligned}
\M{\tau}^\phi_{N_{nn}} &= n_x^2 \M{\tau}_{N_{xx}} +n_x n_y \M{\tau}_{N_{xy}} +n_y^2 \M{\tau}_{N_{yy}} \\
\M{\tau}^\phi_{N_{nt}} &= -n_x n_y \M{\tau}_{N_{xx}} + (n_x^2- n_y^2) \M{\tau}_{N_{xy}} +n_y n_x \M{\tau}_{N_{yy}} \\
\M{\tau}^\phi_{N_{tt}} &= n_y^2 \M{\tau}_{N_{xx}} -n_x n_y \M{\tau}_{N_{xy}} +n_x^2 \M{\tau}_{N_{yy}} \\
\M{\tau}^\phi_{Q_{n}} &= n_x \M{\tau}_{Q_{x}} + n_y \M{\tau}_{Q_{y}} ,\qquad
\M{\tau}^\phi_{Q_{t}} = -n_y \M{\tau}_{Q_{x}} + n_x \M{\tau}_{Q_{y}}
\end{aligned}
\end{equation}
Obviously, the same transformation can be used for moments as well as for tensions.
\end{remark}

\begin{figure}[htb]\centering
   \begin{tikzpicture}[scale=0.6]
\begin{semilogyaxis}[
		xlabel={Iteration $i$},xtick={1,2,3},
		ylabel={\Large $\langle \| \M{\tau}^{Q_x}_{i} - \M{\tau}^{Q_x}_{i-1} \|\rangle/\langle \| \M{\tau}^{Q_x}_{0}  \|\rangle $},
		ylabel near ticks,tick scale binop=\times,
		scaled y ticks=true,]
		\addplot[color=red,mark=x,line width=1.5] table[x index=0,y index=1]{figures/relev/conv.txt};
\end{semilogyaxis}
\end{tikzpicture}
\caption{The norm of the correction to the Saint-Venant stress fields after each iteration}
\label{fig.conv.tau}
\end{figure}
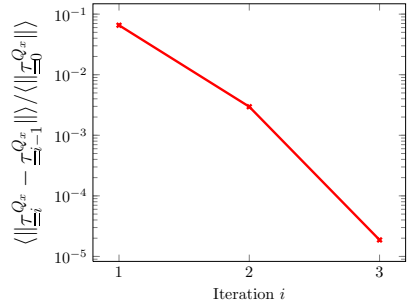

\begin{remark}
The quality of the recovered Saint-Venant stress fields can be improved iteratively: the calculated stress fields are applied as traction boundary conditions to the cell problem, resulting in a problem with even more localized edge effects and a better stress distribution through the thickness in the center of the domain. Figure \ref{fig.conv.tau} shows that the stress field converges very rapidly because each correction is an order of magnitude smaller than the previous one.
\end{remark}

\subsection{Recovery of the Saint-Venant warping vectors}\label{ssec:warping}
Now let us try to define relevant information concerning the displacement field. Since the plate's kinematics (rigid section) is unrealistic and would lead to significant errors \cite{o.allix_plate_2009}, we propose to correct this kinematics using warping vectors (\emph{i.e.} 3D displacements which are allowed to vary through the thickness) deduced from the Saint-Venant auxiliary problems.

From the eight numerical experiments, we obtain the stress states $(\M{\sigma}_j)$ and the displacement profiles $(\V{u}_j)$. By inverting the matrix of the generalized forces of the plate $[\forceP_j^i]$, one obtains the family $(\M{\tau}_i)$ of the uncoupled stresses. The same matrix can be used to define the displacements $\V{U}_i$ associated with pure Saint-Venant loads:
\begin{equation}
\begin{bmatrix} \V{U}_1 &\hdots &\V{U}_8 \end{bmatrix} = \begin{bmatrix} \V{u}_1 &\hdots &\V{u}_8\end{bmatrix} \left[\forceP_j^i\right]^{-1}
\end{equation}
Each of these displacements consists of a plate part and a warping part. In order to separate the two contributions, we assume that the plate displacements produce the same work as the 3D displacements in the Saint-Venant tractions.
{At each point of interface $\inte$}, one has:
\begin{equation}
\begin{aligned}
\dofP_i^k &= \langle \M{\tau}_k:(\V{U}_i\otimes \V{e}_x)\rangle \qquad (i\in\llbracket 1,8\rrbracket,\ k\in\llbracket 1,5\rrbracket)\\
\V{w}_i &= \V{U}_i - \sum_{k=1}^5\shapeP_k \dofP_i^k
\end{aligned}
\end{equation}

\begin{figure}[htb]\centering
\subfloat[][Shear warping, normal component]{
   \begin{tikzpicture}[scale=0.7]
\begin{axis}[
		xlabel={\Large $w^{Q_x}_x$},
		ylabel={thickness z},
		ylabel near ticks,tick scale binop=\times,
		scaled y ticks=true,]
		\addplot[color=red,mark=x,line width=1.5] table[x index=1,y index=0]{figures/relev/ortho_warp_QX.txt};
		\addlegendentry{orthotropic}
		\addplot[color=blue,line width=1.5] table[x index=1,y index=0]{figures/relev/iso_warp_QX.txt};
		\addlegendentry{isotropic}
\end{axis}
\end{tikzpicture}
\label{fig.mu.DX}
}
\subfloat[][Shear warping, tangential component]{
\begin{tikzpicture}[scale=0.7]
\begin{axis}[
		xlabel={\Large $w^{Q_x}_y$},
		ylabel={thickness z},
		ylabel near ticks,tick scale binop=\times,
		scaled y ticks=true,]
		\addplot[color=red,mark=x,line width=1.5] table[x index=2,y index=0]{figures/relev/ortho_warp_QX.txt};
		\addlegendentry{orthotropic}
		\addplot[color=blue,line width=1.5] table[x index=2,y index=0]{figures/relev/iso_warp_QX.txt};
		\addlegendentry{isotropic}
\end{axis}
\end{tikzpicture}
\label{fig.mu.DY}
}
\caption[]{The extracted warping functions for isotropic and orthotropic materials}\label{fig.mu}
\end{figure}
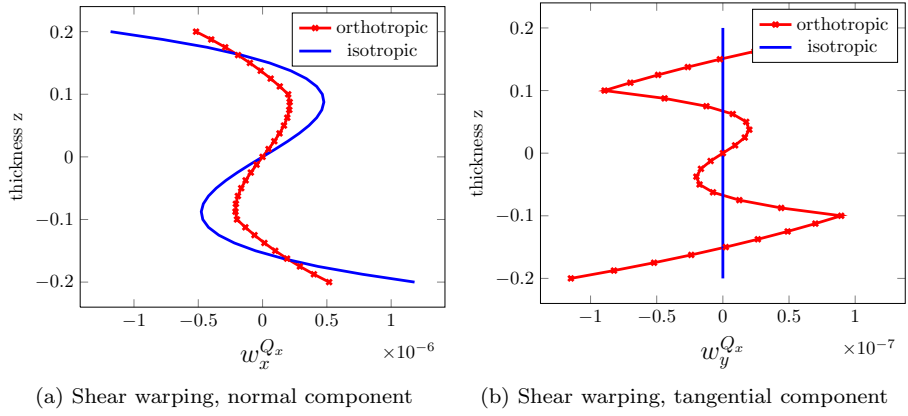

Figure~\ref{fig.mu} presents some of the warping functions extracted from the numerical 3D recovery. The result is in the form of classical zigzags \cite{carrera02} with problem-specific shapes.

\begin{remark}\label{rem:stv_rel}
The extraction of the five kinematic plate components is carried out using only the first five Saint-Venant stress tensors $\V{e}_x$: $(\M{\tau}_k\cdot\V{e}_x)_{k=1:5}$ in the normal direction. Nevertheless, this definition is sufficient for the warping to be well-defined (apart from negligible terms) and orthogonal to all eight Saint-Venant stress tensors for any normal direction. Indeed, from \cite{o.allix_plate_2009}, we know that for Saint-Venant problems the plate displacement corrected by warping in direction $\V{e}_z$ results in 3D stresses which are $O(h/L)$-accurate{, $2h$ being the thickness of the plate}; additional warping in the plane leads to $O(h^2/L^2)$-accuracy. This leads to the following orders of magnitude (which, in practice, turn out to be very pessimistic because the numerical measurements are much smaller):

\begin{equation}\label{eq:stv_rel}
\left.\begin{aligned}
 \langle(\M{\tau}_1\cdot \V{e}_x-  \M{\tau}_2\cdot \V{e}_y)\cdot \V{U}_i \rangle =  \langle(\M{\tau}_2\cdot \V{e}_x- \M{\tau}_6\cdot \V{e}_y)\cdot \V{U}_i \rangle&= o(h^2/L^2) (\Mt{N}_j:\Mt{\gamma}_i)\\
 \langle(\M{\tau}_3\cdot \V{e}_x -  \M{\tau}_4\cdot \V{e}_y)\rangle = \langle(\M{\tau}_4\cdot \V{e}_x - \M{\tau}_7\cdot \V{e}_y) \cdot \V{U}_i\rangle &= o(h^2/L^2)(\Mt{M}_k:\Mt{\chi}_i)\\
 \qquad  \langle(\M{\tau}_5\cdot \V{e}_x -  \M{\tau}_8\cdot \V{e}_y ) \cdot \V{U}_i \rangle&= o(h/L)(\V{Q}_l\cdot\Vt{\epsilon}_{z,i})\end{aligned}\right\}\begin{aligned}\forall i\in\llbracket 1,8\rrbracket\\
\forall j\in\{1,2,6\}\\
\forall k\in\{3,4,7\}\\
\forall l\in\{5,8\}\\
\end{aligned}
\end{equation}
These relations imply that the warping vectors are independent of the axes chosen and produce zero work in all Saint-Venant stress fields, even those that are rotated ($\M{\tau}_k^\phi$) because of linearity.
\end{remark}

\section{Hybrid models}\label{sec:models}
\begin{figure}[ht]\centering\includegraphics[width=.8\textwidth]{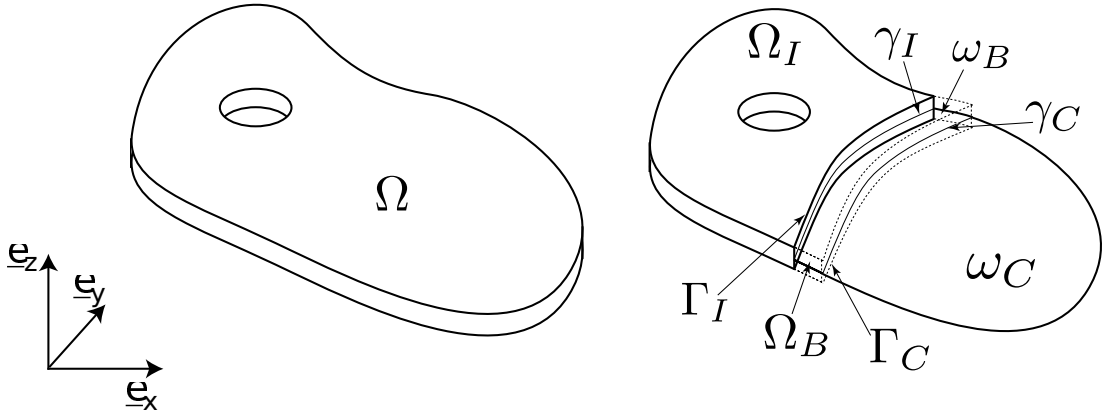}\caption{The reference 3D model and the hybrid model}\label{fig:modeles}\end{figure}

In this section, we focus on the definition of hybrid models for thin structures using plate theories wherever appropriate, but keeping 3D models when necessary (see Figure~\ref{fig:modeles}).
In our context of using commercial software, only less-than-optimal plate theories are available and, thus, a direct plate/3D connection would be error-prone. Besides, the analysis of an industrial structure often begins with a global plate model, so we want to take advantage of this model to the greatest possible extent and supplement it with local 3D modeling only where necessary.

By convention, we use lowercase for the plate domain and for the interfaces ($\omega$, $\gamma$) and uppercase for their 3D counterparts ($\Omega=\omega\times[-h,h]$,  $\Gamma=\gamma\times[-h,h]$). We distinguish three zones:
the zone of interest $I$, the buffer zone $B$ and the complementary zone $C$ (which exists only in the plate model).

\begin{figure}[ht]
\centering
\includegraphics[width=.75\textwidth]{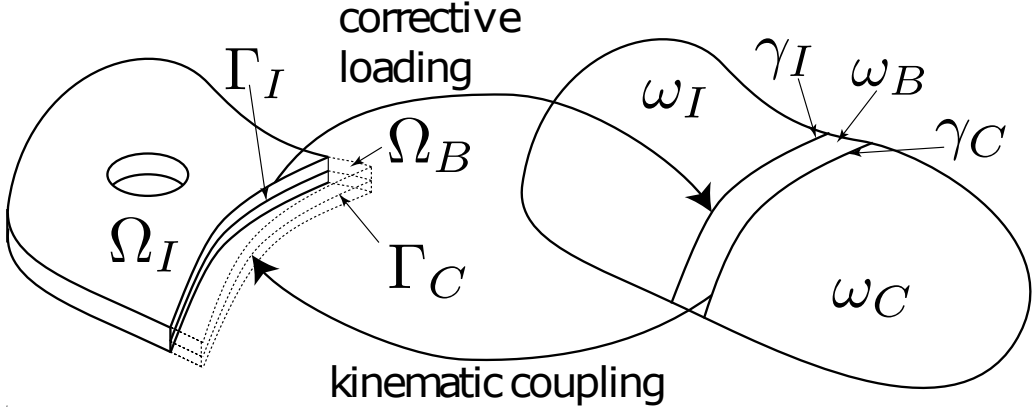}\caption{Nonintrusively coupled models}\label{fig:iteratif}
\end{figure}

Then, the hybrid 3D plate model is defined as the limit of an iterative process, illustrated in Figure~\ref{fig:iteratif}, which alternates between global plate calculations (domain $\domain=\domI\cup\domB\cup\domC$) and local 3D calculations (domain $\DomI\cup\DomB$).


First, a global plate calculation is used to define the boundary conditions for the local 3D problems along $\intC$ (the zooming step). As explained in Section~\ref{ssec:zoom}, the plate-to-3D recovery uses the stress and warping distributions through the thickness deduced from the cell problems of Section~\ref{sec:recovery}.

After the local calculation, the imbalance, in a plate sense, between the plate model and the 3D model is evaluated along the inner line $\intI$ (Section~\ref{ssec:residual}). Then, if necessary, a global plate calculation under a corrective loading applied along $\intI$ is carried out (Section~\ref{ssec:global_correction}).

Beside the dimensionality transition problem, our approach differs from the algorithm of \cite{GENDRE.2009.1} in that we propose to distinguish interface $\intC$ (which serves to prescribe the boundary conditions of the zone of interest) from interface $\intI$ (which is used to apply the global correction) and introduce a buffer zone $\DomB$ in-between.
For the descent steps described in Sections \ref{sec:pure_displacement} and \ref{sssec:zoomLag}, it is possible to choose $\intC=\intI$, but spurious edge effects may occur near $\IntC$. The buffer zone $\DomB$ is used to let the evanescent components activated by the boundary conditions become negligible before entering the zone of interest $\DomI$ itself. In a sense, the plate model and the 3D model are made equivalent in the buffer zone, which makes this method different from other approaches based on a volume connection zone, such as the transition element method \cite{garusi_hybrid_2002,schiermeier_interface_1997} or the Arlequin method \cite{dhia_multiscale_1998}. Besides, the implementation of the method is simpler thanks to its compatibility with nonintrusive approaches.

\subsection{The zooming step}\label{ssec:zoom}
In this section, we assume that a plate problem has been solved and that all its kinematic and static quantities are known along line $\intC$. Let $(\forceP_i)_{i\in\llbracket 1,8\rrbracket}$ and $(\dofP_k)_{k\in\llbracket 1,5\rrbracket}$ denote respectively the 8 static plate components and the 5 displacement components at the interface, and let $s$ be the curvilinear abscissa of interface  $\intC$. Thus, the position of a point of $\IntC$ is given by $(s,z)$ and the classical 3D plate displacement is $\V{v}^P(s,z)=\sum_k\shapeP_k(z)\dofP_k(s)$.

The objective of the zooming step is to define relevant boundary conditions along $\IntC=\intC\times[-h,h]$ in order to solve the 3D problem in the extended zone of interest $\DomI\cup\DomB$ (see Figure~\ref{fig:modeles}). In order to do that, we use the stress and warping distributions $(\M{\tau}_i)$ and $(\V{w}_i)$s calculated for the cell problems of Section~\ref{sec:recovery}.

\subsubsection{Pure traction descent}
At each point of interface $\IntC$, the following traction distribution is prescribed:
\begin{equation}
\M{\sigma}(s,z)\cdot\V{n}(s)=\sum_{i=1}^8 \forceP_i(s)\M{\tau}_i(z) \cdot \V{n}(s) \quad\text{ along }\IntC
\end{equation}
Then, the 3D problem becomes:
\begin{equation}
\begin{aligned}
&\text{Find }\V{u}\in \mathcal{U}(\DomI) \text{ such that }\forall \V{u}^*\in\mathcal{U}_0(\DomI)  \\
&\int\limits_{\DomI} \left(\mathbb{H}:\M{\epsilon}(\V{u})\right):\M{\epsilon}(\V{u}^*)dx = \int\limits_{\DomI} \V{f}\cdot\V{u}^*dx +  \int\limits_{\partial_N\DomI} \V{g}\cdot\V{u}^*dS +\int\limits_{\IntC}\sum_{i=1}^8 \forceP_i(s)(\M{\tau}_i(z) \cdot \V{n}(s))\cdot\V{u}^* dsdz
\end{aligned}
\end{equation}
Of course, this technique works only for domains of interest possessing sufficient Dirichlet conditions. Therefore, it will be evaluated only briefly.

\subsubsection{Lagrangian descent}\label{sssec:zoomLag}
The Lagrangian technique is a kinematic coupling method based on an energy equivalence which is similar to that proposed in  \cite{mccune_mixed-dimensional_2000,shim_mixed_2002} for isotropic plates with known analytic solutions in order to calculate the distribution of Saint-Venant stresses $(\M{\tau}^\phi_i)$. Here, since we are considering stacks of orthotropic plies, no such formula is available and, therefore, we propose to use the numerically calculated $(\M{\tau}^\phi_i)$.

The zooming step consists in prescribing the 3D displacement at the interface which develops the same amount of work as the given plate displacement with the Saint-Venant stresses $(\M{\tau}^\phi_i)$.
\begin{equation}
\begin{aligned}
&\text{Find }\V{u}\in \mathcal{U}(\DomI),\ (\lambda_k) \in\Lambda(\intC)^5 \,\text{such that }\forall \V{u}^*\in\mathcal{U}_0(\DomI) \\
&
\begin{aligned}
\int_{\DomI} \mathbb{H}:\M{\epsilon}(u):\M{\epsilon}(\V{u}^*)dx &= \int_{\DomI}  \V{f}\cdot\V{u}^*dx +  \int_{\partial_f\DomI} \V{g}\cdot\V{u}^*dS \\&+ \sum_{k=1}^5\int_{\intC} \lambda_k(s)\int_z \V{u}^*\cdot\M{\tau}^\phi_k(z)\cdot\V{n}(s)\, dz \, ds\\
& + \sum_{i=6}^8 \int_{\IntC}\forceP^{\phi}_i(s)(\M{\tau}^\phi_i(z) \cdot \V{n}(s))\cdot \V{u}^* dzds \end{aligned}
\\
&  \int_{\intC} \lambda^*_k(s) \int_z (\V{u}-\V{v}^P)\cdot\M{\tau}^{\phi}_k(z)\cdot\V{n}(s)\, dz \, ds=0,\qquad \forall k\in\llbracket 1,5\rrbracket,\forall\lambda^*_k \in\Lambda(\intC)
\end{aligned}
\end{equation}
where $\Lambda(\intC)$ is the space of the Lagrange multipliers $(\lambda_k)$.

In other words, the 5 kinematic degrees of freedom are made {consistent} with their 3D counterparts. The plate part of the 3D displacement, like warping in Section~\ref{ssec:warping}, is defined using averages weighted by $(\M{\tau}^\phi_k\cdot\V{n})$.
The Lagrange multipliers $(\lambda_k)$ are the magnitudes of the normal forces and moments due to the reaction of the 3D domain to the prescribed displacements:
\begin{equation}
k\in\llbracket 1,5\rrbracket,\ \int_z \M{\sigma}(s,z)\cdot\V{n}(s)\cdot\shapeP_k^\phi(z) dz = \lambda_k(s) \int_z\M{\tau}^\phi_k(z)\cdot\V{n}(s)\cdot\shapeP_k^\phi(z) dz = \lambda_k(s)
\end{equation}
The information concerning the last three stress components $\llbracket 6,8\rrbracket$, which involve only tangential action, is given by a classical Neumann condition.

This Lagrangian coupling has the advantage of being applicable to zones of interest with no Dirichlet condition. It also makes the plate's stress components in the 3D zone of interest directly available, which is useful in certain configurations of the global correction step (typically when no buffer zone is used). Unfortunately, this coupling requires the local orientation of the interface to be taken into account, which makes it hard to implement for curved domains. (In the above equations, $\phi$ depends implicitly on $s$.) Also, it is difficult to choose a proper space for multiplier $\Lambda(\intC)$ in the case of nonconforming meshes: the LBB condition must be satisfied (see Section~\ref{ssec:issues_lag}).

\subsubsection{Pure displacement descent}\label{sec:pure_displacement}
Pure displacement descent is a very easy coupling to set up. The principle is to transmit as much information as possible using only Dirichlet conditions at the interface. This eliminates the potential instability problems of Lagrangian coupling and makes the coupling independent of the local orientation of the interface. The prescribed 3D displacement at the interface is given by:
\begin{equation}
\V{u}(s,z)= \sum_{k=1:5} \shapeP_k(z) \dofP_k(s) + \sum_{i=1:8} \V {w}_i(z) \forceP_i(s), \ s\in\intC,\ z\in[-h,h]
\end{equation}

In other words, the classical plate displacements are enriched by warping vectors weighted by the generalized force components of the associated plate.

\subsection{Evaluation of the residual}\label{ssec:residual}
Now let us assume that a 3D problem has been solved through a Lagrangian or pure displacement descent step. This leaves us with a plate solution in $\domain=\domI\cup\domB\cup\domC$ and a 3D solution in $\DomI\cup\DomB$ with kinematic continuity along $\intC$ (in a sense which depends on the descent chosen).
The residual $(\mathbb{L}^\phi_k)$ is defined as the imbalance (in a plate sense) between the plate and the 3D domains along line~$\intI$:
\begin{equation}
\begin{aligned}
\forceP^{\phi,P}_k(s) &: \text{ generalized plate forces along }\intI \text{ viewed as the border of }\domC\cup\domB \\
\forceP^{\phi,3D}_k(s) &= \langle \M{\sigma}:(\V{n}\otimes\shapeP_k) \rangle  \text{: 3D generalized forces along }\intI \text{ viewed as the border of }\domI \\
\mathbb{L}^\phi_k &= \forceP^{\phi,3D}_k+\forceP^{\phi,P}_k \quad\text{ along } \intI,\  k\in\llbracket1,5\rrbracket
\end{aligned}
\end{equation}
{The residual is evaluated in a plate sense because the intention is to reuse it as a corrective load applied to the plate model. A 3D residual could be defined through a recovery of the plate forces, but, besides being unsuitable for the global/local strategy,
it would be difficult to evaluate because an error of order $\left(\frac{h}{L}\right)$ is expected in the $\sigma_{xz}$ and  $\sigma_{yz}$ components and an error of order $\left(\frac{h^2}{L^2}\right)$ is expected in the $\sigma_{zz}$ component.}
{The fact that the 3D quantities have the correct orders of magnitude is interesting information which contributes to the validation of the hybrid model.}

\subsection{The global correction step}\label{ssec:global_correction}
The global correction step consists in solving the global plate problem after having updated the loading by residual $(R_k)$ along interface $\intI$. The updating takes the following form:
\begin{equation}
\sum_{k=1}^5 \int_{\intI}  \mathbb{L}^\phi_k\dofP_k^{\phi,*} ds
\end{equation}
where $\dofP_k^{\phi,*}$ are the degrees of freedom of the plate test displacement field in the coordinate system of $\intI$.

The iterations converge as long as the plate model of the zone of interest $\domI$ is stiffer than the 3D model $\DomI$ (classical fixed point). Otherwise, relaxation is necessary. Such an algorithm also converges for nonlinear monotonic problems \cite{GENDRE.2009.1}.

In the linear case, the convergence can be improved, as is classically done in Schwarz methods, by using Krylov acceleration \cite{widlund_2005}. In this case, the full plate model ($\domain$) can be interpreted as a preconditioner to the hybrid model described in the following section.

In cases where the zooming step has sufficient Dirichlet conditions (\emph{i.e.} the Lagrangian and pure displacement cases), the bufferless limit case $\intI=\intC$ is acceptable. In the Neumann case, a nonzero buffer is necessary. In any case, overlapping is interesting with all methods because it reduces the small edge effects generated by the interface loading of the zooming step.

\subsection{Characterization of the limits}\label{ssec:modellimite}
At the convergence of the iterative process, the residual is zero and the hybrid model can be interpreted as follows:
\begin{itemize}
\item In the zone of interest $\DomI$, it is a pure 3D model,
\item In the complementary zone $\domC$, it is a pure plate model,
\item Along interface $\intI$, the stresses are in equilibrium,
\item Along interface $\intC$, the plate's kinematic continuity is satisfied (along with some other relations which depend on the zooming technique chosen),
\item In the buffer zone, the 3D problem in $\DomB$ is equivalent to the plate problem in $\domB$ in the following sense:
\begin{equation}
\begin{aligned}
\int_{\domB}\left(\langle\M{\sigma}:(\shapeP_k^\phi\otimes \V{e}_\phi)\rangle - \forceP^{P,\phi}_k \right)\dofP_k^{*}dx &= 0  \qquad \forall k\in \llbracket 1,5\rrbracket,\ \forall\phi,\ \forall \dofP_k^*\\
\langle\M{\sigma}\cdot(\shapeP_k^\phi\otimes \V{n})\rangle - \forceP^{P,\phi}_k &= 0 \qquad \text{on } \intI, \ \forall k\in \llbracket 1,5\rrbracket\\
\langle \V{u} \cdot(\M{\tau}^\phi_k \cdot \V{n}) \rangle - \dofP^{P,\phi}_{k}&=0 \qquad \text{on } \intC, \forall k\in \llbracket 1,5\rrbracket
\end{aligned}
\end{equation}
In other words, the plate model obtained by integrating the 3D model through the thickness in ${\DomB}$ must be identical to the plate model in ${\domB}$.
\end{itemize}
One should note that in all cases 3D displacements orthogonal to the Saint-Venant tractions are likely to develop in $\DomB$, but these should be evanescent and should have negligible effect on $\DomI$.

\subsection{Implementation into a finite element code}\label{sec:implementation}
\subsubsection{Preprocessing of the Saint-Venant stresses and warping vectors}
This preprocessor requires an auxiliary 3D problem to be built with the same stacking sequence as the original problem, but a simpler geometry. Minimal Dirichlet conditions are prescribed, except for shear-dominated experiments ($j=5$ and $j=8$).

An important property is that because the plate's kinematics (spanned by functions $(\shapeP_k)$) is linear through the thickness, it is represented exactly by the 3D FE approximation. This means that the calculation of the plate's generalized stresses can be performed directly on the nodal reactions (which are classical outputs of FE software). Let $(\mathbf{F}^l)$ be the nodal reactions through the thickness calculated on a face with normal $\V{e}_x$. If $z^l$ is the height of the $l^{th}$ node), one has:
\begin{equation}
\begin{aligned}
N_{xx} &= \langle \sigma_{xx} \rangle = \sum_l F_x^l  &
N_{xy} &= \langle \sigma_{xy} \rangle = \sum_l F_y^l  & \\
M_{xx} &= \langle z\sigma_{xx} \rangle = \sum_l z^l F_y^l  &
M_{xy} &= \langle z\sigma_{xy} \rangle = -\sum_l z^l F_x^l  & \\
Q_{x} &= \langle \sigma_{xz} \rangle = \sum_l  F_z^l  & &&
\end{aligned}
\end{equation}

\subsubsection{Issues related to the Lagrangian approach}\label{ssec:issues_lag}
The implementation of the Lagrangian approach requires great care in choosing the finite element approximation spaces.

In the case of matching meshes and interpolations (typically 8-node hexahedra on the 3D side and linear plate elements on the plate side), each plate node faces a column of 3D nodes across the thickness. A simple choice consists in using five multiple-point constraints (MPCs) per plate node and corresponding 3D nodes.

The case of nonmatching meshes or interpolations raises the classical issues of mixed approaches: one must check the LBB-stability of the scheme with a well-chosen discretization of the Lagrange multiplier field along the interface.

For example, for a linear plate element and a bilinear 3D element (20-node hexahedron), one method is to use two plate elements facing one column of 3D elements (so each plate node faces a column of 3D nodes), then use a constant Lagrange multiplier for each column of 3D elements.


\section{Numerical study of the various zooming steps}\label{sec:study}

\begin{figure}[ht]
	\centering
    \includegraphics[width=.8\linewidth]{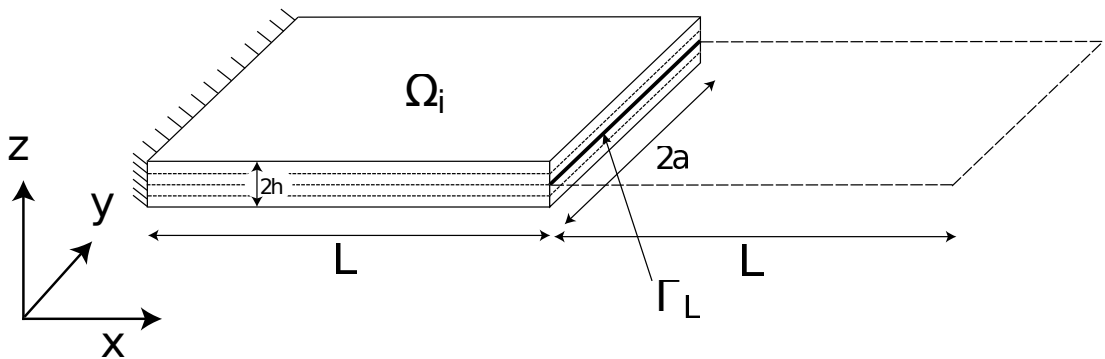}
    \caption{The cantilever plate}
    \label{fig:local_model}
\end{figure}

In order to illustrate the descent step, let us consider the cantilever plate of Figure~ \ref{fig:local_model}, with $L= a = 5$ mm, $h=0.2$ mm {and thickness $2h$}. The plate consists of four orthotropic plies $(-45/45)_{s}$ with the following mechanical properties:
\begin{align*}
E_{L}&=25\,\text{MPa}\\
E_{T}=E_{N}&=1\,\text{MPa}\\
G_{LT}=G_{LN}&=0.5\,\text{MPa}\\
G_{TN}&=0.2\,\text{MPa}\\
\nu_{LT}=\nu_{TN}=\nu_{LN}& = 0.25
\end{align*}

The prescribed displacement at the edge of the plate is $\V{u}_{d}\cdot \V{e}_z=2h$.
The zone of interest is the clamped side because this is where 3D localization takes place. The 3D model was meshed with quadratic hexahedral HEXA20 elements (20 nodes). { A regular mesh was used with 4 elements per ply in the $z$-direction and 20 elements in the $x$ and $y$ directions}. The global model defined over the entire plate was obtained using Dhatt and Batoz homogenization \cite{batoz_moderation_1990} and a discrete shear formulation.
The reference is a full 3D model over the whole domain. The calculations were carried out using Python scripts interfaced with Code\_Aster.

\begin{figure}[ht]\centering
	\subfloat[Extraction along the x-axis]{
\includegraphics[width=.45\textwidth]{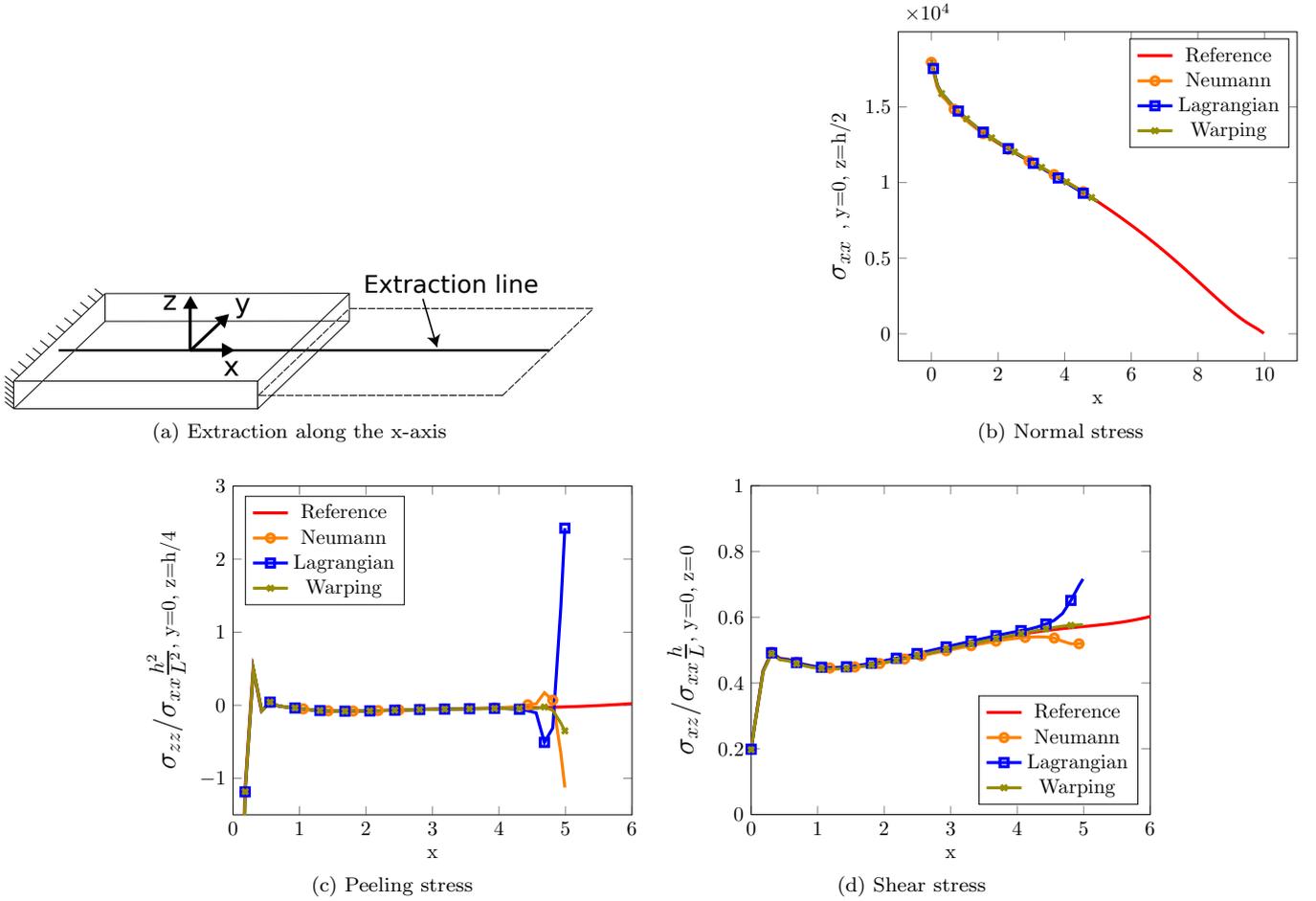}
	}\hfill
\subfloat[][Normal stress]{
	    \begin{tikzpicture}[scale=0.8]
		\begin{axis}[
		xlabel={x},
		ylabel={{\Large $\sigma_{xx}$ },  y=0, z=h/2},
		ylabel near ticks,tick scale binop=\times,
		scaled y ticks=true,]
		\addplot[color=red,line width=1.5] table[x index=0,y index=1]{figures/ortho_45_-45/SIXX_ref_long.txt};
		\addlegendentry{Reference}		
		\addplot[color=orange,mark=o,line width=1.5,mark repeat=6,mark phase=0] table[x index=0,y index=1]{figures/ortho_45_-45/SIXX_force_long.txt};
		\addlegendentry{Neumann}
		\addplot[color=blue,mark=square,line width=1.5,mark repeat=6,mark phase=2] table[x index=0,y index=1]{figures/ortho_45_-45/SIXX_lag_long.txt};
		\addlegendentry{Lagrangian}
		\addplot[color=olive,mark=x,line width=1.5,mark repeat=6,mark phase=4] table[x index=0,y index=1]{figures/ortho_45_-45/SIXX_warp_long.txt};
		\addlegendentry{Warping}
		\end{axis}
		\end{tikzpicture}
		\label{fig:resu_ortho_sigma1}
}

\subfloat[][Peeling stress]{
    	\begin{tikzpicture}[scale=0.8]
		\begin{axis}[
		xlabel={x},
		ylabel={{\Large $\sigma_{zz} / \sigma_{xx}\frac{h^2}{L^2}$}, y=0, z=h/4},
		ylabel near ticks,tick scale binop=\times,
		scaled y ticks=true,
		xmin=0,		xmax=6,
		ymin=-1.5,		ymax=3,
		legend style={legend pos = north west}]
		\addplot[color=red,line width=1.5] table[x index=0,y index=2]{figures/ortho_45_-45/newSIZZ_ref_long.txt};
		\addlegendentry{Reference}
		\addplot[color=orange,mark=o,line width=1.5,mark repeat=3] table[x index=0,y index=2]{figures/ortho_45_-45/newSIZZ_force_long.txt};
		\addlegendentry{Neumann}
		\addplot[color=blue,mark=square,line width=1.5,mark repeat=3] table[x index=0,y index=2]{figures/ortho_45_-45/newSIZZ_lag_long.txt};
		\addlegendentry{Lagrangian}
		\addplot[color=olive,mark=x,line width=1.5,mark repeat=3] table[x index=0,y index=2]{figures/ortho_45_-45/newSIZZ_warp_long.txt};
		\addlegendentry{Warping}
		\end{axis}
		\end{tikzpicture}
      	\label{fig:resu_ortho_sigma2}
}
\subfloat[Shear stress]{
  \begin{tikzpicture}[scale=0.8]
		\begin{axis}[
		xlabel={x},
		ylabel={{\Large $\sigma_{xz} / \sigma_{xx}\frac{h}{L}$}, y=0, z=0},
		ylabel near ticks,tick scale binop=\times,
		scaled y ticks=true,
		xmin=0,		xmax=6,
		ymin=0,		ymax=1,
		legend style={legend pos = south east}]
		\addplot[color=red,line width=1.5] table[x index=0,y index=2]{figures/ortho_45_-45/newSIXZ_ref_long.txt};
		\addlegendentry{Reference}
		\addplot[color=orange,mark=o,line width=1.5,mark repeat=3] table[x index=0,y index=2]{figures/ortho_45_-45/newSIXZ_force_long.txt};
		\addlegendentry{Neumann}
		\addplot[color=blue,mark=square,line width=1.5,mark repeat=3] table[x index=0,y index=2]{figures/ortho_45_-45/newSIXZ_lag_long.txt};
		\addlegendentry{Lagrangian}
		\addplot[color=olive,mark=x,line width=1.5,mark repeat=3] table[x index=0,y index=2]{figures/ortho_45_-45/newSIXZ_warp_long.txt};
		\addlegendentry{Warping}
	\end{axis}
	\end{tikzpicture}
	\label{fig:resu_ortho_sigma3}}
\caption[]{Evolution of stresses along the x-axis} \label{fig:resu_ortho_sig}
\end{figure}

Since the zone of interest is clamped, it was possible to use all three zooming techniques. Figure~\ref{fig:resu_ortho_sig} shows the evolution of some components of the 3D stress along the $x$ axis. The reference is defined in $\Domain$ ($x\in[0,10]$) while the 3D parts of the hybrid models are defined in zone $x\in[0,5]$.
In all three methods, as expected theoretically \cite{o.allix_plate_2009}, $\sigma_{xz}$ has approximately the same order of magnitude as $\sigma_{xx}(\frac{h}{L})$ and $\sigma_{zz}$ has approximately the same order of magnitude as $\sigma_{xx}(\frac{h^{2}}{L^{2}})$.
One can observe that in all cases a perturbation appears near the interface ($x=5$). However, since the generalized forces transmitted at the interfaces are correct, this local effect decreases rapidly and vanishes within a $2h$ band. The 3D solutions are close to the reference solution. Altogether, the perturbations caused by the warping technique seem smaller than those caused by the Neumann technique, which, in turn, are smaller than with Lagrangian coupling. The presence of these undesirable stress concentrations is what triggered the extension of the zone of interest by a 3D buffer zone in order to dampen the artificial edge effects.

\begin{figure}[ht]	
\subfloat[Extraction along the y-axis]{
\includegraphics[width=.45\textwidth]{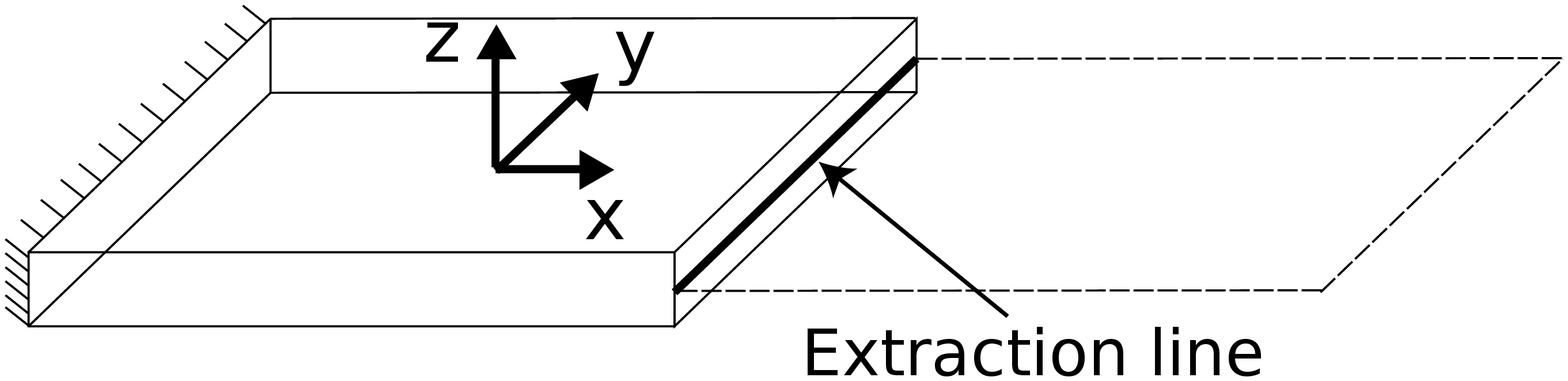}
	}
\subfloat[Bending moment]{
	\begin{tikzpicture}[scale=0.8]
		\begin{axis}[
		xlabel={y},
		ylabel={\Large $M_{xx}$},
		ylabel near ticks,tick scale binop=\times,
		scaled y ticks=true,]
		\addplot[color=red,line width=1.5] table[x=COOR_Y,y=MXX]{figures/ortho_45_-45/plaque.txt};
		\addlegendentry{Plate}		
		\addplot[color=orange,mark=o,line width=1.5,mark repeat=3] table[x=COOR_Y,y=MXX]{figures/ortho_45_-45/3D_force.txt};
		\addlegendentry{Neumann}
		\addplot[color=blue,mark=square,line width=1.5,mark repeat=3] table[x=COOR_Y,y=MXX]{figures/ortho_45_-45/3D_lag.txt};
		\addlegendentry{Lagrangian}
		\addplot[color=olive,mark=x,line width=1.5,mark repeat=3] table[x=COOR_Y,y=MXX]{figures/ortho_45_-45/3D_warp.txt};
		\addlegendentry{Warping}
		\end{axis}
		\end{tikzpicture}
	\label{fig:resu_ortho_EG1}
}

\subfloat[Torsion moment]{
	\begin{tikzpicture}[scale=0.8]
		\begin{axis}[
		xlabel={y},
		ylabel={\Large $M_{xy}$},
		ylabel near ticks,tick scale binop=\times,
		scaled y ticks=true,
		legend style={legend pos = north west}]
		\addplot[color=red,line width=1.5] table[x=COOR_Y,y=MXY]{figures/ortho_45_-45/plaque.txt};
		\addlegendentry{Plate}		
		\addplot[color=orange,mark=o,line width=1.5,mark repeat=3] table[x=COOR_Y,y=MXY]{figures/ortho_45_-45/3D_force.txt};
		\addlegendentry{Neumann}
		\addplot[color=blue,mark=square,line width=1.5,mark repeat=3] table[x=COOR_Y,y=MXY]{figures/ortho_45_-45/3D_lag.txt};
		\addlegendentry{Lagrangian}
		\addplot[color=olive,mark=x,line width=1.5,mark repeat=3] table[x=COOR_Y,y=MXY]{figures/ortho_45_-45/3D_warp.txt};
		\addlegendentry{Warping}
		\end{axis}
		\end{tikzpicture}
	\label{fig:resu_ortho_EG2}
}
\subfloat[Normal shear]{
	\begin{tikzpicture}[scale=0.8]
	\begin{axis}[
		xlabel={y},
		ylabel={\Large $Q_{x}$},
		ylabel near ticks,tick scale binop=\times,
		scaled y ticks=true,
		legend style={legend pos = south east}]
		\addplot[color=red,line width=1.5] table[x=COOR_Y,y=QX]{figures/ortho_45_-45/plaque.txt};
		\addlegendentry{Plate}		
		\addplot[color=orange,mark=o,line width=1.5,mark repeat=3] table[x=COOR_Y,y=QX]{figures/ortho_45_-45/3D_force.txt};
		\addlegendentry{Neumann}
		\addplot[color=blue,mark=square,line width=1.5,mark repeat=3] table[x=COOR_Y,y=QX]{figures/ortho_45_-45/3D_lag.txt};
		\addlegendentry{Lagrangian}
		\addplot[color=olive,mark=x,line width=1.5,mark repeat=3] table[x=COOR_Y,y=QX]{figures/ortho_45_-45/3D_warp.txt};
		\addlegendentry{Warping}
	\end{axis}
	\end{tikzpicture}
	\label{fig:resu_ortho_EG3}}
\caption{The generalized forces along the interface}\label{fig:resu_ortho_forcY}
\end{figure}

Figure~\ref{fig:resu_ortho_forcY} shows the evolution of the generalized stresses along the interface for the initial plate calculation and for the 3D calculations. By construction, Neumann zooming produces the same generalized stresses as the plate calculation. Conversely, the other zooming techniques lead to slightly different results, which justifies the iterations between the local and global models. As expected, the differences between the plate model and the zooms are particularly significant near the edges.


\section{Evaluation of the hybrid model based on warping}\label{sec:assess}
In this section, we study the application of the iterative coupling technique in order to achieve convergence toward the hybrid model. Based on previous experiments, we chose the warping-based descent algorithm (see \ref{sec:pure_displacement}) and used a nonzero buffer zone (see \ref{ssec:global_correction}).

\begin{figure}[ht]
	\centering
    \includegraphics[width=.99\linewidth]{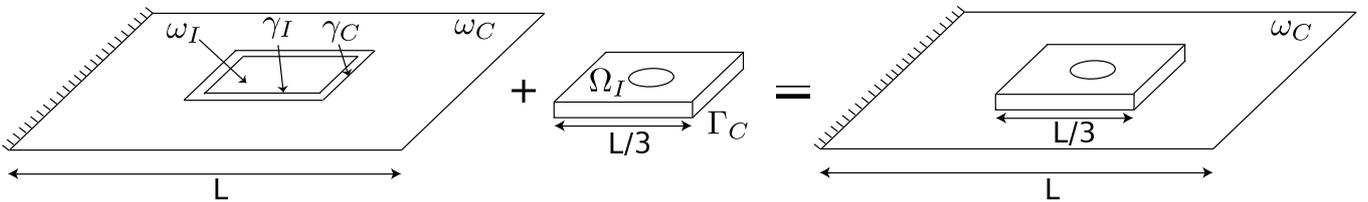}
    \caption{The coupling of the two models (left) defines the hybrid model (right)}
    \label{fig:hole_plate}
\end{figure}

The objective was to model the holed composite plate shown in Figure~\ref{fig:hole_plate}. The local model consisted in a 3D representation of the hole and its vicinity with enough extent to capture the 3D edge effect. The global model was a plate, in which we could choose to represent (or not) the hole in the zone $\domI$ corresponding to the 3D model. We used the same orthotropic layers as in the previous section. The dimensions were $L= a = 60$ mm and $ h=1$ mm. The radius of the hole was $r=2$ mm. The loading was a transverse displacement $u_{d}=0.45$ along the right side.

The zone of interest $\domI$ was a square of side $20$ mm and the width of the buffer zone $\domB$ was $2$ mm. {It was meshed with 24 elements per side and 16 elements in the thickness}.

\subsection{Convergence toward the hybrid model}
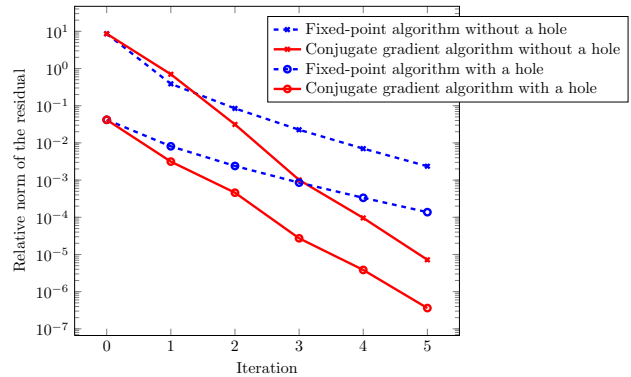
\begin{figure}[ht!]
\centering
\begin{tikzpicture}[scale=0.6]
	\begin{semilogyaxis}[
		scale only axis,		
		xlabel={Iteration},
		ylabel=Relative norm of the residual, ylabel near ticks,tick scale binop=\times,
		yticklabel style={/pgf/number format/fixed, /pgf/number format/precision=5},
		legend style={legend cell align=left, legend pos = north east,anchor=north}]
			\addplot[color=blue,mark=x,line width=1.5,dashed,mark options={solid}] table[x = iteration,y = classic]{figures/conv.txt};
		\addlegendentry{Fixed-point algorithm without a hole}
			\addplot[color=red,mark=x,line width=1.5] table[x = iteration,y = GC]{figures/conv.txt};
		\addlegendentry{Conjugate gradient algorithm without a hole}
e					\addplot[color=blue,mark=o,line width=1.5,dashed,mark options={solid}] table[x = iteration,y = trou_classic]{figures/conv.txt};
		\addlegendentry{Fixed-point algorithm with a hole}
					\addplot[color=red,mark=o,line width=1.5] table[x = iteration,y = trou_GC]{figures/conv.txt};
		\addlegendentry{Conjugate gradient algorithm with a hole}
	\end{semilogyaxis}
	\end{tikzpicture}
	\caption{Evolution of the relative norm of the residual as a function of the number of iterations}
    \label{fig:iteration}
\end{figure}

Let us define the relative norm $\eta_i$ of the residual $\mathbb{L}$ of Section~\ref{ssec:residual} as follows:
\begin{equation*}
\eta_i ^2 = \frac{\int_{\intI}|(\Mt{N}_i^{3D}+\Mt{N}_i^P) \cdot \V{n}|^2 ds}{\max_{\intI}|\Mt{N}_0^P \cdot \V{n}|^2} + \frac{\int_{\intI}|(\Mt{M}_i^{3D}+\Mt{M}_i^P) \cdot \V{n}|^2 ds}{\max_{\intI}|\Mt{M}_0^P \cdot \V{n}|^2} + \frac{\int_{\intI}|(\Vt{Q}_i^{3D}+\Vt{Q}_i^P) \cdot \V{n}|^2 ds}{\max_{\intI}|\Vt{Q}_0^P \cdot \V{n}|^2}
\end{equation*}
Figure~\ref{fig:iteration} shows the decrease in $\eta_i$ as a function of the number of iterations for the classical fixed point algorithm and for the conjugate gradient algorithm. The complexity of the plate model in the zone of interest $\domI$ is a parameter of the method which influences the convergence, but not the limit. We considered two cases, depending on whether the hole was represented in the plate model $\domI$ or not.

When the hole was included, the plate model gave a better representation of the zone of interest, leading to much less initial imbalance. In both cases, the convergence was fast. In particular, one can observe that two conjugate gradient iterations of the model without a hole led to a similar quality approximation as the classical submodeling technique applied to the holed model (which is more difficult to generate and to mesh).

\subsection{Validation of the hybrid model}

Now let us study the hybrid model obtained after a sufficient number of iterations of the coupling procedure and compare the results of this model with those of a full 3D calculation of $\Domain$ in terms of the reliability of the 3D quantities obtained.

\begin{figure}[ht!]
\centering
\subfloat[][In-plane stress]{
\begin{tikzpicture}[scale=0.57]
	\begin{axis}[
		scale only axis,		
		xlabel={$r$, $z=h$},
		ylabel={\Large  $\sigma_{xx}$}, ylabel near ticks,tick scale binop=\times,xmin=0,xmax=20,
		yticklabel style={/pgf/number format/fixed, /pgf/number format/precision=5},
		title={\Large $\sigma_{xx}$ $45^{\circ}$},
		legend style={legend pos = south east},]
		\addplot[color=blue,line width=1.5] table[x = arc_length,y = SIXX]{figures/plaque_trou_ref/SIXX_45.txt};
		\addlegendentry{Reference}
		\addplot[color=red,mark=x,line width=1.5] table[x = arc_length,y = SIXX]{figures/plaque_trou/SIXX_45_0.txt};
		\addlegendentry{Iteration 0}
		\addplot[color=olive,mark=o,line width=1.5] table[x = arc_length,y = SIXX]{figures/plaque_trou/SIXX_45_conv.txt};
		\addlegendentry{Iteration 5}
	\end{axis}
	\end{tikzpicture}
}
\subfloat[][Shear stress]{
\begin{tikzpicture}[scale=0.57]
	\begin{axis}[
		scale only axis,		
		xlabel={$r$, $z=0$},
		ylabel={\Large $\sigma_{xz}$}, ylabel near ticks,tick scale binop=\times,xmin=0,xmax=20,
		yticklabel style={/pgf/number format/fixed, /pgf/number format/precision=5},
		title={\Large  $\sigma_{xz}$ $45^{\circ}$},
		legend style={legend pos = north east},]
		\addplot[color=blue,line width=1.5] table[x = arc_length,y = SIXZ]{figures/plaque_trou_ref/SIXZ_45.txt};
		\addlegendentry{Reference}
		\addplot[color=red,mark=x,line width=1.5] table[x = arc_length,y = SIXZ]{figures/plaque_trou/SIXZ_45_0.txt};
		\addlegendentry{Iteration 0}
		\addplot[color=olive,mark=o,line width=1.5] table[x = arc_length,y = SIXZ]{figures/plaque_trou/SIXZ_45_conv.txt};
		\addlegendentry{Iteration 5}
	\end{axis}
	\end{tikzpicture}
}

\subfloat[][Moment]{
\begin{tikzpicture}[scale=0.57]
	\begin{axis}[
		scale only axis,		
		xlabel={$r$},
		ylabel={\Large  $M_{xx}$}, ylabel near ticks,tick scale binop=\times,xmin=0,xmax=20,
		yticklabel style={/pgf/number format/fixed, /pgf/number format/precision=5},
		title={\Large $M_{xx}$ $45^{\circ}$},
		legend style={legend pos = south east},]
		\addplot[color=blue,line width=1.5] table[x = R_ref, y = MXX_ref, col sep=semicolon,]{figures/conv/EG_45.csv};
		\addlegendentry{Reference}
		\addplot[color=red,mark=x,line width=1.5] table[x = R_loc, y = MXX_0, col sep=semicolon,]{figures/conv/EG_45.csv};
		\addlegendentry{Iteration 0}
		\addplot[color=olive,mark=o,line width=1.5] table[x = R_loc, y = MXX_conv, col sep=semicolon,]{figures/conv/EG_45.csv};
		\addlegendentry{Iteration 5}
		\addplot[color=black,line width=1.5,dashed] table[x = R_MXX, y = MXX_plaque, col sep=semicolon,]{figures/conv/EG_45.csv};
		\addlegendentry{Plate with a hole}
	\end{axis}
	\end{tikzpicture}
}
\subfloat[][Shear force]{
\begin{tikzpicture}[scale=0.57]
	\begin{axis}[
		scale only axis,		
		xlabel={$r$},
		ylabel={\Large $Q_{x}$}, ylabel near ticks,tick scale binop=\times,xmin=0,xmax=20,
		yticklabel style={/pgf/number format/fixed, /pgf/number format/precision=5},
		title={\Large  $Q_{x}$ $45^{\circ}$},
		legend style={legend pos = north east},]
		\addplot[color=blue,line width=1.5] table[x = R_ref, y = QX_ref, col sep=semicolon,]{figures/conv/EG_45.csv};
		\addlegendentry{Reference}
		\addplot[color=red,mark=x,line width=1.5] table[x = R_loc, y = QX_0, col sep=semicolon,]{figures/conv/EG_45.csv};
		\addlegendentry{Iteration 0}
		\addplot[color=olive,mark=o,line width=1.5] table[x = R_loc, y = QX_conv, col sep=semicolon,]{figures/conv/EG_45.csv};
		\addlegendentry{Iteration 5}
		\addplot[color=black,line width=1.5,dashed] table[x = R_QX, y = QX_plaque, col sep=semicolon,]{figures/conv/EG_45.csv};
		\addlegendentry{Plate with a hole}
	\end{axis}
	\end{tikzpicture}
}

\subfloat[][Peeling stress]{
\begin{tikzpicture}[scale=0.565]
	\begin{axis}[
		scale only axis,		
		xlabel={$r$, $z=h/2$},
		ylabel={\Large $\sigma_{zz}$}, ylabel near ticks,tick scale binop=\times,xmin=0,xmax=20,
		yticklabel style={/pgf/number format/fixed, /pgf/number format/precision=5},
		title={\Large $\sigma_{zz}$ $45^{\circ}$},
		legend style={legend pos = north east},]
		\addplot[color=blue,line width=1.5] table[x = arc_length,y = SIZZ]{figures/plaque_trou_ref/SIZZ_45.txt};
		\addlegendentry{Reference}
		\addplot[color=red,mark=x,line width=1.5] table[x = arc_length,y = SIZZ]{figures/plaque_trou/SIZZ_45_0.txt};
		\addlegendentry{Iteration 0}
		\addplot[color=olive,mark=o,line width=1.5] table[x = arc_length,y = SIZZ]{figures/plaque_trou/SIZZ_45_conv.txt};
		\addlegendentry{Iteration 5}
	\end{axis}
	\end{tikzpicture}
	}
\subfloat[][The stress extraction along a $45^{\circ}$ line]{
	\includegraphics[width=.5\linewidth]{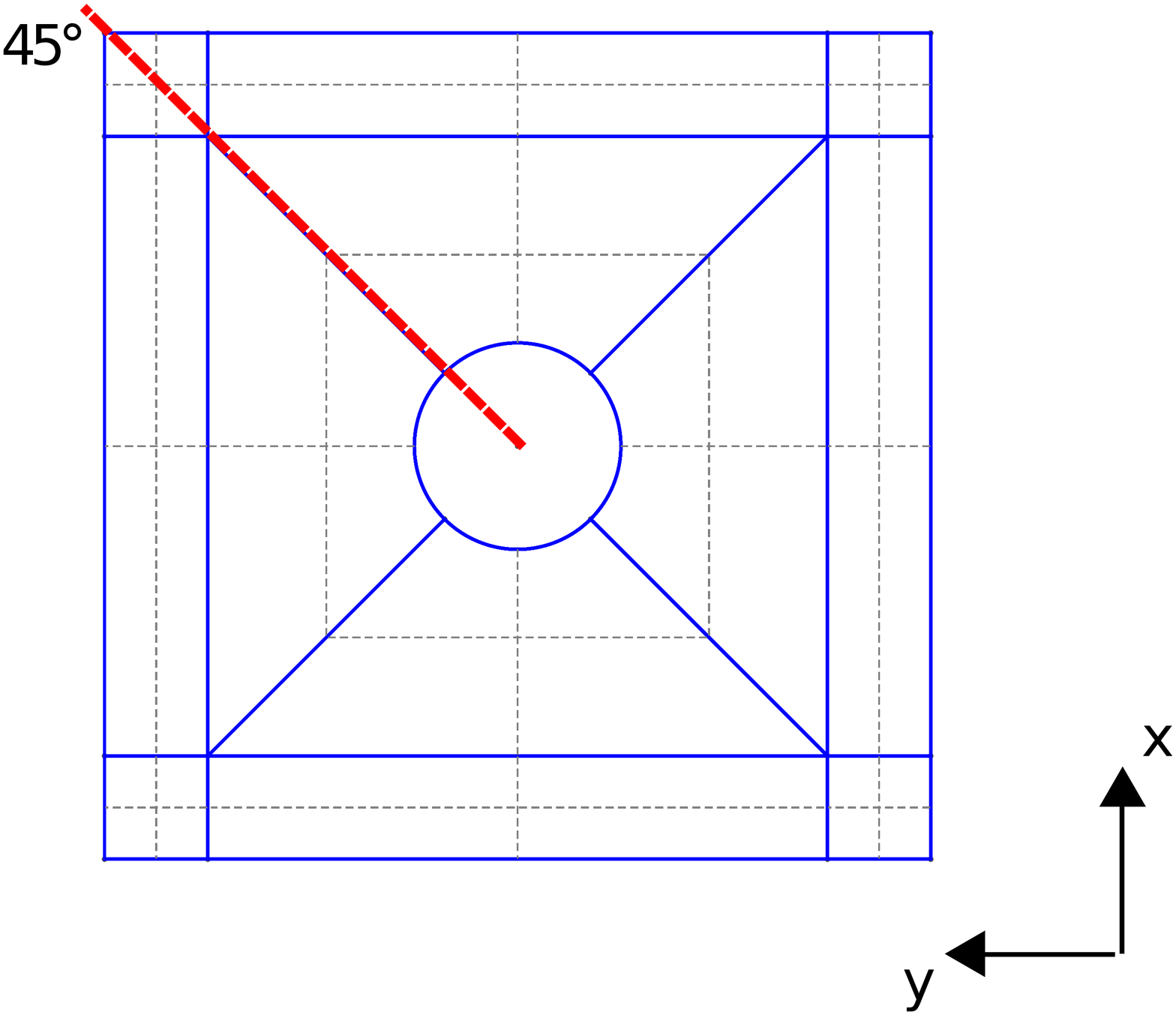}\label{fig:extract_line}
	}
	\caption{Stress and generalized force distributions near the hole}\label{fig:stress_hole}
\end{figure}

The quantities of interest chosen are the stresses in the vicinity of the hole, which would initiate delamination.
Figures~\ref{fig:stress_hole}(a,b,e) show the variations of some stress components along the particular line described in Figure~\ref{fig:extract_line}.
A comparison of the full 3D reference values, the values using submodeling (which is equivalent to the hybrid model with a single descent and no iteration) and the values given by the hybrid model after 5 iterations clearly shows that submodeling is very unreliable while the hybrid model is almost identical to the reference model.

Figures~\ref{fig:stress_hole}(c,d) show the evolution of the generalized forces along the same line. A comparison of the reference 3D calculation, the hybrid model at iterations 0 and 5 and a direct (holed) plate model shows that both the plate and the submodel lead to poor results while the hybrid model comes very close to the reference.

Table~\ref{tab:plaque_conv} quantifies the changes induced by the iterations in the plate quantities in the plate domain: the final values of the deflection, shear stress and moment (\emph{i.e.} the values corresponding to the hybrid model) at a point of $\intI$ ($x=20$ mm, $y=0$ mm) are compared with the initial values obtained with the plate model (with or without a hole). Of course, the model without a hole had significant errors, but even the model with a hole was improved by taking into account 3D phenomena. (In particular, the shear stress was underestimated by more than $5\%$.)

\begin{table}[htbp]
\centering
\pgfplotstableset{create on use/---/.style={create col/set list={error in $w$ (\%), error in $Q_{x}$ (\%), error in $M_{xx}$ (\%)},{string type}}}
\pgfplotstabletypeset[
columns={---,A,B},
col sep=&,
empty cells with={--},
columns/---/.style={string type},
columns/A/.style={column name=with a hole,precision=2},
columns/B/.style={column name=without a hole,precision=2},
every head row/.style={before row=\toprule,after row=\midrule},
every last row/.style={after row=\bottomrule},]
	{
	A&B
-0.0193227&	-1.9667018
5.9926455&	-35.9376882	
1.1771289&	3.2961126
	}
	\caption{The corrections between the first descent and the fifth iteration}\label{tab:plaque_conv}
\end{table}


\section{Conclusion}\label{sec:conclu}
In this paper, we proposed a new strategy for modeling thin composite structures by coupling a global plate model with a refined 3D model in zones of interest. The objective is to combine the efficiency of plate modeling for most of the structure with the added accuracy of 3D calculations where necessary (near edges, holes or defects). The hybrid plate/3D model is obtained after the convergence of a nonintrusive iterative process which enables one to use standard commercial software and industrial plate models without even having to represent the critical zones. The coupling is achieved by means of a preliminary analysis in which numerical Saint-Venant stress and warping fields which are suitable for the stacking sequence and the chosen discretization are recovered automatically. Several coupling methods are possible, but the warping-based technique is the simplest to implement and seems to lead to the least amount of perturbation at the interface.
The convergence can be monitored by looking at the imbalance (in a plate sense) between the 3D domain and the plate domain. In the case of linear 3D domains, it can be accelerated by using a Krylov solver.

Our works in progress concern the application of the method to the simulation of complex assemblies of composite plates. In this case, the global model uses classical simplified 1D connectors to maintain computational efficiency while realistic 3D nonlinear models of bolts, rivets and damageable composites are used in the local models. A particular issue arises from the use of different in-plane discretizations at the interface between the local model and the global model. Two-scale approaches, such as that of \cite{guidault2007two} in the context of the XFEM, have interesting advantages, such as the possibility to transfer average static plate quantities to the 3D model without introducing local forces and spurious edge effects.

\section*{Acknowledgement}
This work was partially funded by the French National Research Agency as part of project ICARE (ANR-12-MONU-0002-04).

\bibliography{article}

\end{document}